\documentclass[]{interact}
\usepackage{latexsym}
\usepackage{anyfontsize} 
\usepackage{tikz} 
\usetikzlibrary{cd} 
\usepackage{color}
\usepackage{hyperref}
\usepackage{url}
\usepackage{breakurl}
\newcommand{\bburl}[1]{\textcolor{blue}{\url{#1}}}
\usepackage[backend=biber, style=alphabetic, maxnames=99]{biblatex} 
\makeatletter
\newcommand{\monthyear}[1]{%
  \def\@monthyear{\uppercase{#1}}}
\newcommand{\volnumber}[1]{%
  \def\@volnumber{\uppercase{#1}}}
\makeatother

\newcommand{\Z}{{\mathbb Z}}

\theoremstyle{plain}
\numberwithin{equation}{section} 
\newtheorem{thm}{Theorem}[section] 
\newtheorem{theorem}[thm]{Theorem}
\newtheorem{lemma}[thm]{Lemma}

\newtheorem{definition}[thm]{Definition}
\newtheorem{proposition}[thm]{Proposition}
\newtheorem{corollary}[thm]{Corollary}
\newtheorem{remark}[thm]{Remark}

\numberwithin{table}{section} 
\numberwithin{figure}{section}

\begin{document}

\monthyear{August 2026}
\volnumber{Volume, Number}
\setcounter{page}{1}

\title{Random Fibonacci Tilings}

\author{
\name{Yuval Amit\textsuperscript{a} (yuval.amit33@gmail.com), Azmera Gebre\textsuperscript{b} (a.gebre0813@gmail.com), Christopher Housholder\textsuperscript{c} (christopherlancehousholder@gmail.com), Joshua Khan\textsuperscript{d} (joshuakhan1307@gmail.com), Steven J. Miller\textsuperscript{b} (sjm1@williams.edu), Emily Upton\textsuperscript{e} (emilyupton220@gmail.com), Connor Yau\textsuperscript{b}(cy10@williams.edu)
\thanks{This work was supported financially by the National Science Foundation (grant DMS2341623), Williams College Finnerty Fund, the University of Michigan, Missouri State University, and the University of Rochester. We want to thank Williams College and SMALL REU for continuous support.
}}
\affil{\textsuperscript{a} Department of Mathematics,
                University of Maryland, College Park, Maryland, USA
      	\\ \textsuperscript{b} Department of Mathematics, Williams College, Williamstown, Massachusetts, USA 
        \\ \textsuperscript{c} Department of Mathematics,
                Missouri State University, 
                Springfield, Missouri, USA
        \\ \textsuperscript{d} Department of Mathematics, University of Rochester, Rochester, New York, USA
        \\ \textsuperscript{e} Department of Mathematics, University of Rhode Island, Kingston, Rhode Island, USA}
}

\maketitle

{\bf Article type}: research 
\bigskip

\begin{abstract}
The Fibonacci spiral is formed by a square-joining rule: begin with a unit square and attach squares whose side lengths follow the Fibonacci sequence; in the standard construction, choices of direction produce a spiral and a tiling of the plane. We study a randomization of this construction, in which the attachment axis alternates deterministically, while fair coin flips determine direction. For a lattice point we ask when the growing rectangle first covers that point and how the expected covering time depends on its location. Tracking the rectangle’s boundary reduces this question to random subset sums of even- and odd-indexed Fibonacci numbers. We extend the construction to $d$ dimensions, where coordinate subsequences are super-increasing, allowing us to derive exact covering probabilities and formulas for the expectation and variance. For $d\geq2$, the expected covering time grows logarithmically with distance from the initial block, while the variance remains uniformly bounded. We obtain exponential tail bounds for the delay beyond the earliest feasible covering stage. Finally, we extend the Fibonacci spiral to the randomized tiling through a square-by-square construction using circular arcs and quintic connectors. This construction uses one-step lookahead, admits a continuously differentiable parametrization, and recovers the Fibonacci spiral for the standard tiling.
\end{abstract}


\section{Introduction}\label{sec:introduction}

\begin{figure}[h]
    \centering
    \includegraphics[width=0.47\textwidth]{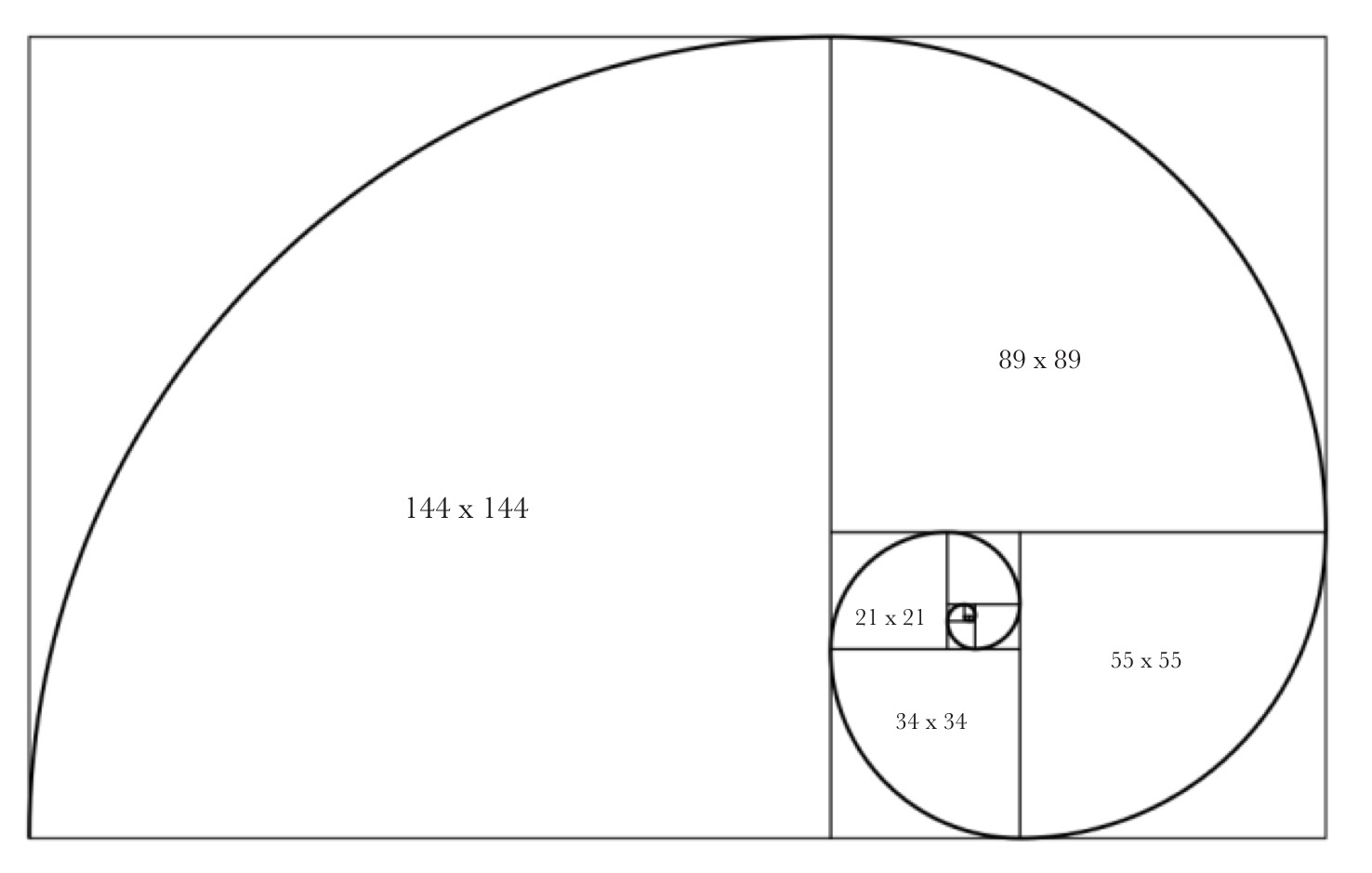}
    \caption{Fibonacci Spiral Construction}
    \label{fig:fibonacci-spiral}
\end{figure}

The Fibonacci sequence is one of the most recognizable and fascinating sequences in mathematics. It's defined recursively by \begin{equation}
    F_0 \ = \ F_1 \ = \ 1,
    \qquad F_n \ = \ F_{n-1}+F_{n-2}\quad(n\geq2),
\end{equation} and arises in many and diverse settings, from pure mathematics to natural biological structures \cite{Kos01}. One of its more interesting properties is that if we create a square tile for each Fibonacci number, then these $F_n \times F_n$ cover the plane. We place the tiles one at a time in the pattern right - up - left - down. In doing so we draw a quarter-circle from from the previous through the next tile, as shown in Figure ~\ref{fig:fibonacci-spiral} 

In  ~\cite[Section~2]{EM}, the authors explore a game where players begin with a rectangle, rather than the unit square, and add squares in a similar fashion to the Fibonacci Spiral, and observe the geometric patterns that emerge. 

Instead of changing the original starting shape from a $1 \times 1$ square to a $1 \times n$ rectangle, we consider a new generalization. We study a random version of the square-seed construction. Starting from the fixed unit square, we alternate the attachment axis and use independent fair coin flips to choose the side on which each new square is added. Thus the pattern is now (right or left) - (up or down) - (right or left) - (up or down) and so on, where each alternative happens with probability one-half. Our main question is when the growing rectangle first covers a given lattice point. After solving that, we extend the spiral construction to these random placements, as the random alternations complicates creating such a curve (we'll see we won't know how to add the curve in a tile unless we know the positioning of the next).

In the covering problem, our main simplification is to track the rectangle’s boundary. The total boundary displacement at a given time is a random subset of fixed weights. Furthermore, the covering time of each coordinate is independent, so the total covering time is found just as their maximum. In a cyclic extension to $d$ dimensions, the weights are subsequences of a generalized Fibonacci recurrence. The superincreasing property of these generalized sequences orders the subset sums and yields exact covering probabilities, expectations, and variances.

Theorem~\ref{thm:uniform-fluctuations} states that for every target with \(I(a)\neq\varnothing\) and every integer \(k\geq0\),
\begin{equation}
\label{eq:uniform-delay-tail}
\mathbb P(T\ >\ B(a)+dk)\ \leq\ \min\{1,|I(a)|2^{-k}\}\ \leq\ \min\{1,d2^{-k}\}.
\end{equation}
For every fixed real \(q\geq1\),
\begin{equation}
\label{eq:uniform-delay-moments}
\mathbb E[|T-B(a)|^q]\ =\ O_{d,q}(1),
\end{equation}
uniformly over the target.

For each fixed $d\geq2$, Corollary~\ref{cor:covering-asymptotics} gives
\begin{equation}
    \mathbb{E}[T(a)]
      \ = \ \frac{\log\delta_\infty(a)}{\log\rho_d}+O_d(1),
    \qquad
    \operatorname{Var}(T(a)) \ = \ O_d(1)
    \label{eq:intro-covering-asymptotics}
\end{equation}
as $\delta_\infty(a)\to\infty$. Here $T(a)$ counts moves after the
initial cube $[0,1]^d$, $\rho_d>1$ is the positive root of
$x^d-x^{d-1}-1$, and
\begin{equation}
    \delta_r(a) \ :=\ \max\{a_r-1,-a_r,0\},
    \qquad
    \delta_\infty(a) \ := \ \max_{0\leq r<d}\delta_r(a).
    \label{eq:intro-target-distance}
\end{equation}
where $a_r$ is the $r$\textsuperscript{th} coordinate of the target point $a$. Thus the mean grows logarithmically with distance while the variance remains bounded. We also obtain exponential tails for the delay beyond the earliest feasible covering stage.


\section{Random Tiling Construction}
\label{sec:random-tiling}
\subsection{The planar construction}

We use the square-seed rectangle game of~\cite[Section~2]{EM}, with
the following indexing and placement rule. Let
\begin{equation}
    F_0 \ = \ F_1 \ = \ 1,
    \qquad F_n \ = \ F_{n-1}+F_{n-2}\quad(n\geq2).
    \label{eq:game-fibonacci}
\end{equation}
Our $F_n$ is the Fibonacci number denoted $F_{n+1}$ in
\cite[Definition~2.4]{EM}. 

We work in the Euclidean plane \(\mathbb R^2\). Let $R_n^{(2)}\subset\mathbb R^2$ be the closed rectangle after $n$ moves, starting with
$R_0^{(2)}=[0,1]^2$. The initial square is not counted as a move.
Let $(c_n)_{n\geq1}$ be independent bits with
\begin{equation}
    \mathbb{P}(c_n=0) \ = \ \mathbb{P}(c_n=1)=\frac12.
    \label{eq:game-placement-bits}
\end{equation}
At move $n$, append an $F_n\times F_n$ square along a full side of
the current rectangle. For odd $n$, place it on the right if $c_n=0$
and on the left if $c_n=1$. For even $n$, place it above if $c_n=0$
and below if $c_n=1$. Thus $0$ selects right/up and $1$ selects
left/down.

The rectangle update in~\cite[Proposition~2.2 and Remark~2.5]{EM}
gives side lengths
\begin{equation}
    F_n\quad\text{and}\quad F_{n+1}
    \label{eq:planar-side-lengths}
\end{equation}
after $n$ moves, up to their orientation. The random choices affect
the position of the rectangle, while its dimensions remain deterministic.

\begin{figure}[ht]
    \centering
    \includegraphics[width=0.6\linewidth]{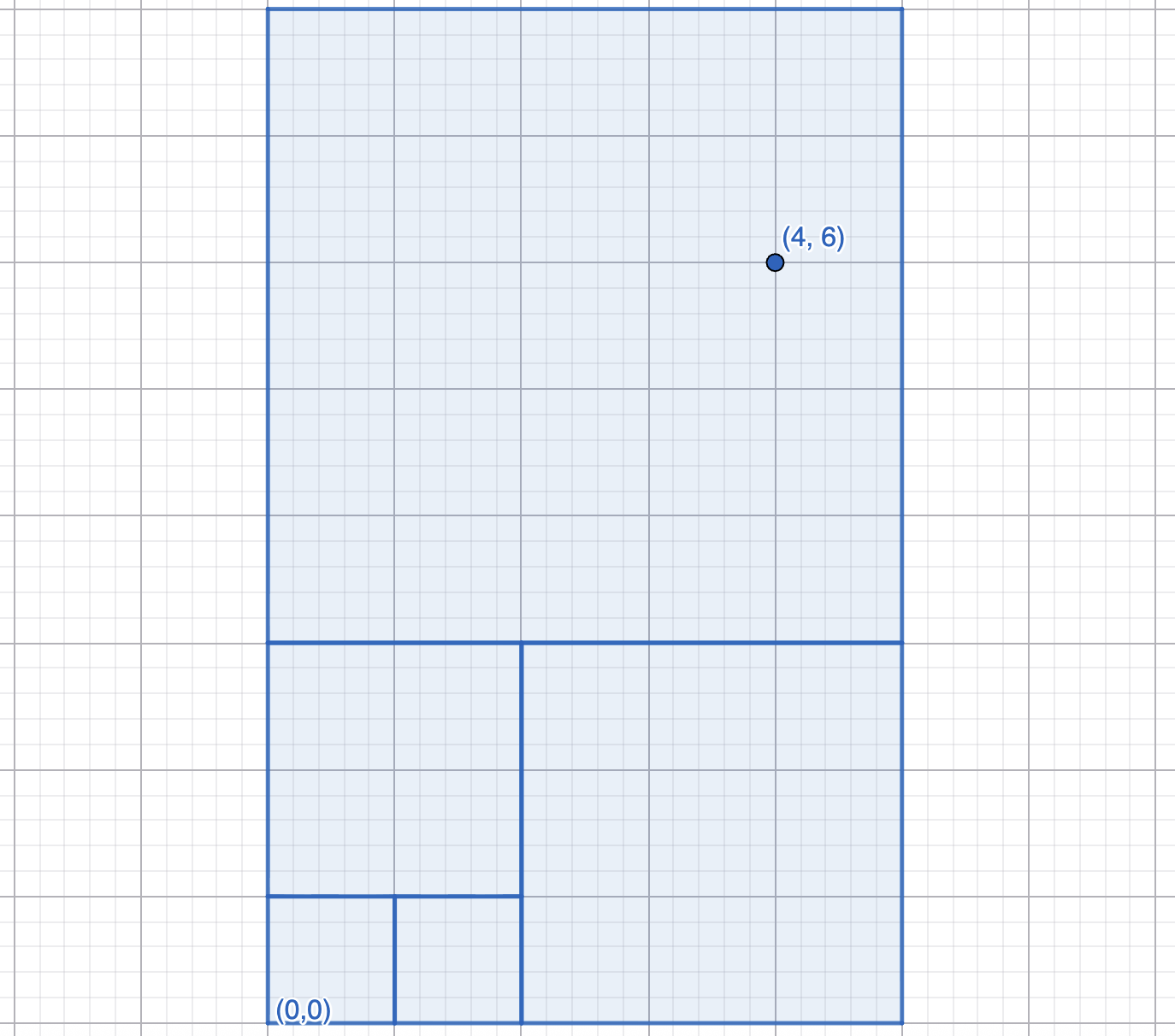}
    \caption{Attaining $(4,6)$ in four moves.}
    \label{ex1}
\end{figure}

\begin{figure}[ht]
    \centering
    \includegraphics[width=0.6\linewidth]{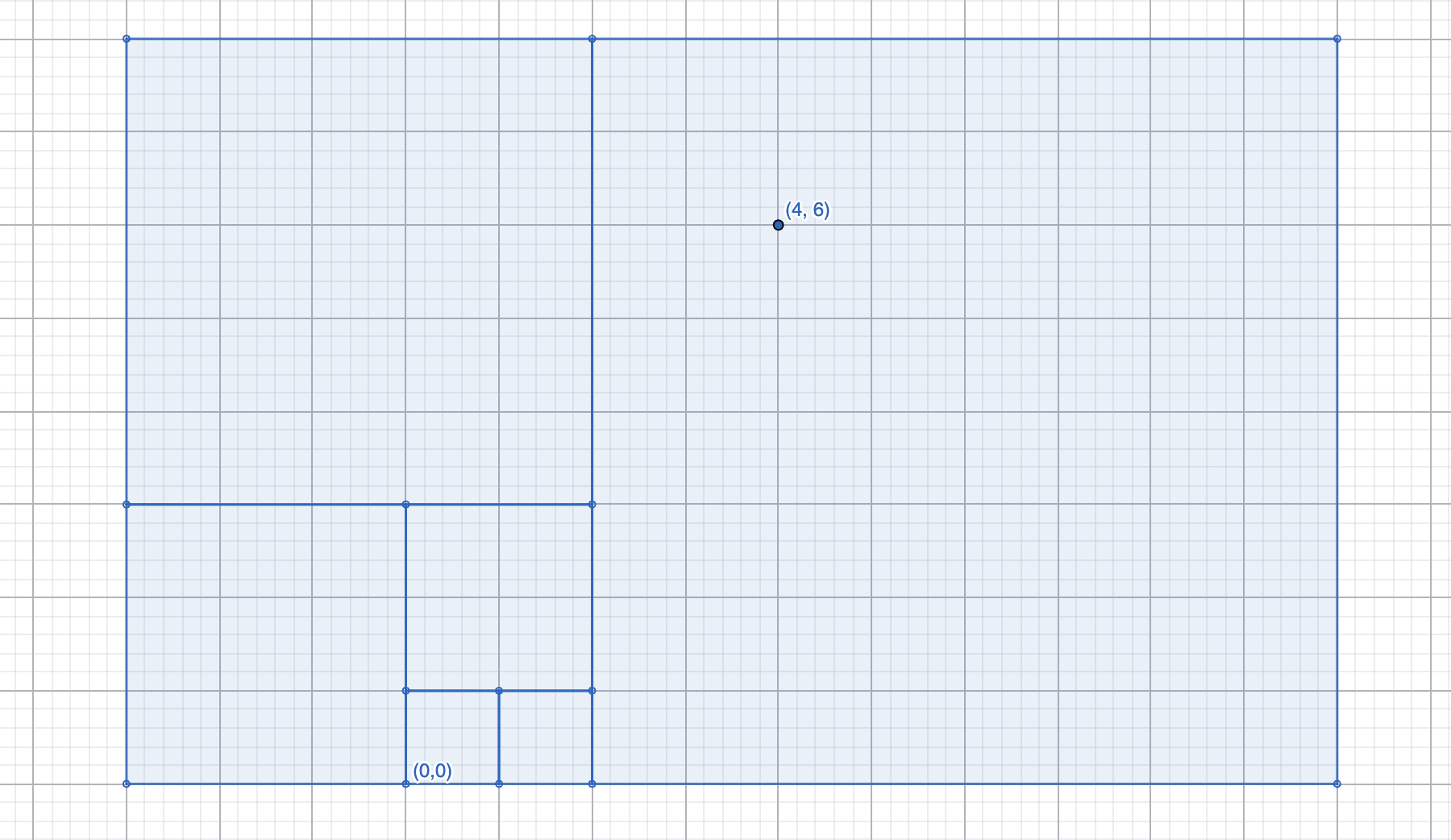}
    \caption{Attaining $(4,6)$ in five moves.}
    \label{ex2}
\end{figure}

For $(a,b)\in\mathbb{Z}^2$, its covering time is
\[
    T(a,b)=\inf\{n\geq0:(a,b)\in R_n^{(2)}\}.
\]
Boundary points count as covered, and points in the initial square have
covering time zero. The placement words $(0,0,0,0)$ and $(0,0,1,0,0)$
in Figures~\ref{ex1} and~\ref{ex2} illustrate coverage of $(4,6)$ after four and
five moves, respectively. Section~\ref{sec:subset-sums} develops the subset sums needed
to compute such covering times.

The construction gives rise to two principal questions.  First, for a lattice point $(a,b)\in\Z^2$, how many moves are required, on average, before the growing rectangle covers $(a,b)$?  Second, how can the Fibonacci spiral from the standard deterministic tiling be extended to this random setting?

\subsection{The higher-dimensional construction}

The planar process has a natural cyclic extension to dimension $d\geq 2$. Begin with the unit hypercube $[0,1]^d$. At move $k\geq1$, update coordinate $r\equiv k-1\pmod d$, with $0\leq r<d$. Use the independent fair placement bits $(c_k)_{k\geq1}$, where $0$ selects the positive face and $1$ selects the negative face.  At successive moves, choose the coordinate axes in a fixed cyclic order.  Along the selected axis, a fair coin determines whether the new block is attached to the positive or negative face. The new block spans the entire exposed $(d-1)$-dimensional face, and its thickness is equal to the most recently created side length.  Equivalently, at each move the side length that has gone one complete coordinate cycle without being updated is replaced by its sum with the most recently created side length.

This geometric rule motivates the following recurrence. The sequence is the specialization $u(n;d-1,1)$ of Harris and Styles \cite{HS}; its recurrence and spectral properties are classical.

\begin{definition}[Order-$d$ Fibonacci sequence]
\leavevmode\label{def:order-d}
For an integer $d\geq 2$, the order-$d$ Fibonacci sequence
$(F_n^{(d)})_{n\geq 0}$ is defined by
\begin{equation}
\label{eq:order-d-initial-values}
F_0^{(d)}\ =\ F_1^{(d)}\ =\ \cdots\ =\ F_{d-1}^{(d)}\ =\ 1
\end{equation}
and
\begin{equation}
\label{eq:order-d-recurrence}
F_n^{(d)}\ =\ F_{n-1}^{(d)}+F_{n-d}^{(d)},
\qquad n\ \geq\ d.
\end{equation}
\end{definition}

\begin{proposition}
\leavevmode\label{prop:side-lengths}
After $k\geq 0$ moves, the side lengths of the bounding hyperrectangle are
\begin{equation}
\label{eq:cuboid-side-lengths}
F_k^{(d)},\ F_{k+1}^{(d)},\ \ldots\ ,\ F_{k+d-1}^{(d)},
\end{equation}
up to a cyclic permutation of the coordinate axes.

For $k\geq 1$, immediately before move $k$ the selected coordinate has length $F_{k-1}^{(d)}$, while the remaining face has side lengths
\begin{equation}
\label{eq:exposed-face-dimensions}
F_k^{(d)},\ F_{k+1}^{(d)},\ \ldots\ ,\ F_{k+d-2}^{(d)}.
\end{equation}
The appended block spans this face and has thickness $F_{k+d-2}^{(d)}$. Consequently, the selected side length becomes
\begin{equation}
\label{eq:updated-side-length}
F_{k-1}^{(d)}+F_{k+d-2}^{(d)}\ =\ F_{k+d-1}^{(d)}.
\end{equation}
\end{proposition}

\begin{proof}
Before any moves are made, all $d$ side lengths are equal to $1$, so the claim holds for $k=0$.  Now proceed by induction and suppose it holds immediately before move $k$.  The cyclic axis rule selects the coordinate with length $F_{k-1}^{(d)}$.  The appended block has thickness $F_{k+d-2}^{(d)}$ and spans the full opposite face, so the other $d-1$ side lengths remain unchanged.  The selected side becomes
\begin{equation}
\label{eq:side-length-induction}
F_{k-1}^{(d)}+F_{k+d-2}^{(d)}\ =\ F_{k+d-1}^{(d)}
\end{equation}
by the recurrence.  Thus the new collection of side lengths is
\begin{equation}
\label{eq:induction-side-vector}
F_k^{(d)},\ F_{k+1}^{(d)},\ \ldots\ ,\ F_{k+d-1}^{(d)},
\end{equation}
which completes the induction.
\end{proof}

For $d=2$, the exposed face has length $F_k^{(2)}$ and the appended block has the same thickness, so it is an $F_k^{(2)}\times F_k^{(2)}$ square.  Hence the construction reduces exactly to the planar Fibonacci tiling.  For $d=3$, the recurrence
\begin{equation}
\label{eq:narayana-recurrence}
F_n^{(3)}\ =\ F_{n-1}^{(3)}+F_{n-3}^{(3)}
\end{equation}
is the Narayana recurrence \cite{OEIS930}, and the appended objects are generally rectangular cuboids rather than cubes.

\begin{remark}
\leavevmode\label{rem:scalar-analogue}
For our purposes, in the one-dimensional analogue we do not directly substitute $d=1$ into the higher-dimensional construction. Instead, we study the analogous scalar covering problem while retaining the Fibonacci sequence from the planar model. Our one-dimensional analysis is a probabilistic reduction used to understand a single coordinate, whereas the geometric construction itself is defined for $d\geq 2$.
\end{remark}

We next record the growth properties needed in the covering-time analysis. Many properties of this sequence are intuitively borrowed from the Fibonacci sequence itself.

\begin{proposition}
\leavevmode\label{prop:coordinate-superincreasing}
For all integers \(n,k\geq0\),
\begin{equation}
\label{eq:coordinate-doubling}
F_{n+kd}^{(d)}\ \geq\ 2^kF_n^{(d)}.
\end{equation}
Moreover, for \(k\geq1\),
\begin{equation}
\label{eq:superincreasing-chain}
F_{n+(k-1)d}^{(d)}\ \leq\ \sum_{j=0}^{k-1}F_{n+jd}^{(d)}\ <\ 2F_{n+(k-1)d}^{(d)}\ \leq\ F_{n+kd}^{(d)}.
\end{equation}
Thus every subsequence \((F_{n+jd}^{(d)})_{j\geq0}\) with fixed
starting index \(n\) is superincreasing.
\end{proposition}

\begin{proof}
The initial values and recurrence imply by induction that all terms are
positive, that the sequence is nondecreasing, and that it is strictly
increasing from index \(d\) onward. Hence
\begin{equation}
\label{eq:one-cycle-doubling}
F_{n+d}^{(d)}\ =\ F_{n+d-1}^{(d)}+F_n^{(d)}\ \geq\ 2F_n^{(d)}.
\end{equation}
Iteration proves the first assertion. For \(0\leq j<k\), it also gives
\begin{equation}
\label{eq:earlier-weight-bound}
F_{n+jd}^{(d)}\ \leq\ 2^{-(k-1-j)}F_{n+(k-1)d}^{(d)}.
\end{equation}
Summing these inequalities yields
\begin{equation}
\label{eq:superincreasing-proof}
\sum_{j=0}^{k-1}F_{n+jd}^{(d)}\ \leq\ F_{n+(k-1)d}^{(d)}\sum_{i=0}^{k-1}2^{-i}\ <\ 2F_{n+(k-1)d}^{(d)}\ \leq\ F_{n+kd}^{(d)}.
\end{equation}
The lower bound follows by retaining the largest summand.
\end{proof}

Proposition~\ref{prop:coordinate-superincreasing} is what we actually need for studying covering time. This property allows us to impose an ordering on the subsets of these modular subsequences by comparing the most significant elements.

We now define the dominant root of the recurrence, which provides information about the growth rate of the sequence. 

\begin{definition}
\leavevmode\label{def:dominant-root}
The characteristic polynomial of the order-$d$ Fibonacci recurrence is
\begin{equation}
\label{eq:characteristic-polynomial}
\chi_d(x)\ =\ x^d-x^{d-1}-1.
\end{equation}
Its unique positive real root is denoted by $\rho_d$ and is called the dominant root.  Equivalently,
\begin{equation}
\label{eq:dominant-root-identity}
\rho_d^{d-1}(\rho_d-1)\ =\ 1.
\end{equation}
It satisfies
\begin{equation}
\label{eq:dominant-root-range}
1\ <\ \rho_d\ <\ 2.
\end{equation}
For $d=2$, one recovers the golden ratio
\begin{equation}
\label{eq:golden-ratio}
\rho_2\ =\ \varphi\ =\ \frac{1+\sqrt5}{2}.
\end{equation}
\end{definition}

\begin{lemma}
\leavevmode\label{lemma:spectral}
Every root of $\chi_d$ is simple.  Furthermore, if $\zeta$ is a root of $\chi_d$
and $\zeta\neq\rho_d$, then
\begin{equation}
\label{eq:spectral-gap}
|\zeta|\ <\ \rho_d.
\end{equation}

\end{lemma}

\begin{proof}
If $\zeta$ were a repeated root, then $\chi_d(\zeta)=\chi_d'(\zeta)=0$.
Because
\begin{equation}
\label{eq:characteristic-derivative}
\chi_d'(x)\ =\ x^{d-2}\bigl(dx-(d-1)\bigr)
\end{equation}
and $\chi_d(0)=-1$, the only possible common root would be $\zeta=(d-1)/d$.  However,
\begin{equation}
\label{eq:excluded-repeated-root}
\chi_d\!\left(\frac{d-1}{d}\right)\ =\ -\frac1d\left(\frac{d-1}{d}\right)^{d-1}-1\ <\ 0,
\end{equation}
so no repeated root exists.

Now let $\chi_d(\zeta)=0$.  Then
\begin{equation}
\label{eq:root-modulus-identity}
1\ =\ |\zeta|^{d-1}|\zeta-1|.
\end{equation}
Suppose that $|\zeta|\geq\rho_d$.  Since $\rho_d>1$, the reverse triangle inequality and the strict increase of $r^{d-1}(r-1)$ for $r\geq1$ give
\begin{equation}
\label{eq:strict-dominance-proof}
1\ =\ |\zeta|^{d-1}|\zeta-1|\ \geq\ |\zeta|^{d-1}\bigl(|\zeta|-1\bigr)\ \geq\ \rho_d^{d-1}(\rho_d-1)\ =\ 1.
\end{equation}
Equality must therefore hold throughout.  The second equality forces $|\zeta|=\rho_d$, and equality in $|\zeta-1|\geq|\zeta|-1$ forces $\zeta$ to be nonnegative and real. Hence $\zeta=\rho_d$.  Every other root therefore satisfies $|\zeta|<\rho_d$.
\end{proof}

\begin{theorem}
\leavevmode\label{theorem:asymptotic}
Let $\alpha_d=\max\bigl\{|\zeta|:\chi_d(\zeta)=0,\ \zeta\neq\rho_d\bigr\}$.
Then $\alpha_d<\rho_d$, and as $n\to\infty$, $
F_n^{(d)}=A_d\rho_d^n+O\!\left(\alpha_d^n\right),$
where
\begin{equation}
\label{eq:dominant-coefficient}
A_d\ =\ \frac{1}{\rho_d^{-1}\left(1+d\rho_d^{-(d-1)}\right)}\ =\ \frac{\rho_d^d}{\rho_d^{d-1}+d}\ >\ 0.
\end{equation}
In particular,
\begin{equation}
\label{eq:dominant-equivalence}
F_n^{(d)}\ \sim\ A_d\rho_d^n.
\end{equation}

\end{theorem}

\begin{proof}
The recurrence and initial conditions give the generating function
\begin{equation}
\label{eq:generating-function}
\sum_{n\geq0}F_n^{(d)}z^n\ =\ \frac{1}{1-z-z^d}.
\end{equation}
Let $D_d(z) = 1-z-z^d$. The zeros of $D_d(z)$ are the reciprocals of the roots of $\chi_d$. Label the roots of $\chi_d$ as $\zeta_1, \zeta_2,\dots,\zeta_d$, with $\zeta_d = \rho_d$. Thus
\begin{equation}
\label{eq:denominator-factorization}
D_d(z)\ =\ \prod_{j=1}^d(1-\zeta_jz).
\end{equation}

There exist constants $A_1,\dots,A_d$ such that
\begin{equation}
\label{eq:partial-fractions}
\begin{aligned}
\sum_{n\geq0} F_n^{(d)}z^n\ =\ \frac{1}{D_d(z)} &\ =\ \frac{1}{\prod_{j=1}^d(1-\zeta_jz)}\\
    &\ =\ \sum_{j=1}^d\frac{A_j}{1-\zeta_jz}\\
    &\ =\ \sum_{j=1}^dA_j\sum_{n\geq 0} \zeta_j^nz^n.
\end{aligned}
\end{equation}
Therefore, by extracting the coefficient of $z^n$, 
\begin{equation}
\label{eq:binet-expansion}
F_n^{(d)}\ =\ \sum_{j=1}^dA_j\zeta_j^n.
\end{equation}

By Lemma~\ref{lemma:spectral}, $z=\rho_d^{-1}$ is the unique zero of smallest
modulus, and all zeros are simple. Furthermore, every root satisfies 
$|\zeta_j|\leq\alpha_d<\rho_d$ for $1\leq j\leq d-1$, so
\begin{equation}
\label{eq:nondominant-remainder}
\left|\sum_{j=1}^{d-1} A_j\zeta_j^n\right|\ \leq\ 
\left(\sum_{j=1}^{d-1}|A_j|\right)\alpha_d^n\ =\ O(\alpha_d^n).
\end{equation}

It follows that
\begin{equation}
\label{eq:growth-expansion}
F_n^{(d)}\ =\ A_d\rho_d^n+O(\alpha_d^n).
\end{equation}
Dividing by \(A_d\rho_d^n\), we obtain
\begin{equation}
\label{eq:normalized-growth}
\frac{F_n^{(d)}}{A_d\rho_d^n}\ =\ 1+O\left(\left(\frac{\alpha_d}{\rho_d}\right)^n\right)
\longrightarrow1,
\end{equation}
because \(\alpha_d/\rho_d<1\). Thus
\begin{equation}
\label{eq:growth-equivalence-proof}
F_n^{(d)}\ \sim\ A_d\rho_d^n.
\end{equation}

The coefficient of the dominant term is obtained from the simple pole at
$z_d=\rho_d^{-1}$, giving
\begin{equation}
\label{eq:dominant-residue}
A_d\ =\ -\frac{1}{z_d\,\frac{d}{dz}(1-z-z^d)\rvert_{z=z_d}}\ =\ \frac{1}{z_d(1+dz_d^{d-1})}\ =\ \frac{\rho_d^d}{\rho_d^{d-1}+d}.
\end{equation}
This quantity is positive.
\end{proof}

\begin{corollary}
\leavevmode\label{cor:inverse-growth}
As $n\to\infty$,
\begin{equation}
\label{eq:logarithmic-growth}
\log F_n^{(d)}\ =\ n\log\rho_d+\log A_d+o(1)\ =\ n\log\rho_d+O(1).
\end{equation}
Equivalently, if $n_d(A)=\min\{n\geq0:F_n^{(d)}\geq A\}$ then as $A$ tends to infinity,
\begin{equation}
\label{eq:inverse-growth}
n_d(A)\ =\ \frac{\log A}{\log\rho_d}+O(1)\ =\ \log_{\rho_d}A+O(1).
\end{equation}
\end{corollary}

\begin{proof}
By Theorem~\ref{theorem:asymptotic},
\begin{equation}
\label{eq:growth-factorization}
F_n^{(d)}\ =\ A_d\rho_d^n\bigl(1+o(1)\bigr).
\end{equation}
Taking logarithms proves the first assertion.  Solving the resulting estimate
for the threshold index proves the second.
\end{proof}

\begin{corollary}
\leavevmode\label{cor:planar-growth}
For $d=2$ (i.e., the case in $\Z^2$), the order-$d$ sequence is the ordinary Fibonacci sequence and $\rho_2=\varphi=(1+\sqrt5)/2$. Consequently,
\begin{equation}
\label{eq:planar-binet}
F_n^{(2)}\ =\ \frac{\varphi^{n+1}}{\sqrt5}+O\!\left(\varphi^{-n}\right)
\end{equation}
and
\begin{equation}
\label{eq:planar-threshold-growth}
n_2(A)\ =\ \log_\varphi A+O(1).
\end{equation}
\end{corollary}

\subsection{Covering times}

Let $\mathcal{R}_n^{(d)}$ denote the random bounding hyperrectangle after $n$ completed moves.  For a lattice point $a=(a_0,\ldots,a_{d-1})\in\Z^d$, define
\begin{equation}
\label{eq:target-distances-active-set}
\delta_r(a)\ =\ \max\{a_r-1,-a_r,0\},
\qquad
I(a)\ =\ \{r:\delta_r(a)\ >\ 0\}.
\end{equation}
Thus $\delta_r(a)$ is the distance from the appropriate face of $[0,1]^d$ to the target in coordinate $r$, and $I(a)$ is the set of coordinates not already covered initially.  Define the covering time by
\begin{equation}
\label{eq:covering-time-definition}
T(a)\ :=\ \inf\{n\ \geq\ 0:a\in\mathcal{R}_n^{(d)}\}.
\end{equation}
The first problem of the paper is to determine the distribution and moments of $T(a)$, beginning with
\begin{equation}
\label{eq:covering-time-objective}
\mathbb{E}[T(a)].
\end{equation}
Although the side lengths of $\mathcal{R}_n^{(d)}$ are deterministic, its
position depends on the coin-flip word.  Determining whether a given point has
been covered is therefore equivalent to ordering certain weighted subset sums.
The $j$\textsuperscript{th} update of coordinate $r$ occurs at move $r+1+jd$ and has thickness
\begin{equation}
\label{eq:coordinate-weight-definition}
w_{r,j}\ :=\ F_{d-1+r+jd}^{(d)},
\qquad j\ \geq\ 0.
\end{equation}
For each fixed $r$, these weights form a tail of a residue-class subsequence, so Proposition~\ref{prop:coordinate-superincreasing} supplies the key structural fact for the analysis in the next section.

The second problem, the construction and analysis of a random Fibonacci spiral, uses the same sequence of appended blocks but tracks distinguished corners rather than the entire bounding hyperrectangle.  We save its formal introduction for Section~\ref{sec:random-spirals}.

\section{Weighted Subset Sums and Canonical Ordering}
\label{sec:subset-sums}
The side lengths of the random tiling are deterministic, but the position of the bounding hyperrectangle depends on the directions selected by the coin flips. Along each coordinate axis, the distance traveled toward a fixed target is therefore a weighted subset sum. We show that the terms associated with any one coordinate form a superincreasing sequence. It follows that their subset sums are ordered exactly as binary integers, which yields a canonical threshold word and an exact finite-time covering probability.

Throughout this section, let $\left\{F_n^{(d)}\right\}_{n\geq0}$ be the order-$d$ Fibonacci sequence defined by
\begin{equation}
\label{eq:recurrence-recall}
F_0^{(d)}\ =\ \cdots\ =\ F_{d-1}^{(d)}\ =\ 1, \qquad
F_n^{(d)}\ =\ F_{n-1}^{(d)}+F_{n-d}^{(d)} \text{ for } n\ \geq\ d.
\end{equation}

\subsection{Coordinate subsequences and random subsets}

Fix $r\in\{0,1,\ldots,d-1\}$.  The $j$\textsuperscript{th} update of coordinate $r$, indexed from $j=0$, occurs at move $r+1+jd$ and moves the selected face by
\begin{equation}
\label{eq:coordinate-weight-recall}
w_{r,j}\ :=\ F_{d-1+r+jd}^{(d)}.
\end{equation}
For $m\geq0$, define the first $m$ coordinate weights and their cumulative sum by
\begin{equation}
\label{eq:coordinate-weight-set-sum}
\mathcal W_{r,m}\ :=\ \{w_{r,0},\ldots,w_{r,m-1}\},
\qquad
W_{r,m}\ :=\ \sum_{j=0}^{m-1}w_{r,j},
\end{equation}
with $\mathcal W_{r,0}=\varnothing$ and $W_{r,0}=0$.  Proposition~\ref{prop:coordinate-superincreasing} gives
\begin{equation}
\label{eq:coordinate-superincreasing}
w_{r,k}\ >\ W_{r,k}
\qquad(k\ \geq\ 1),
\end{equation}
so every coordinate-weight sequence is superincreasing.

Let $n\geq0$ denote the number of completed moves.  The number of updates of coordinate $r$ among those $n$ moves is
\begin{equation}
\label{eq:update-count}
u_r(n)\ :=\ \#\{0\ \leq\ i\ <\ n:i\equiv r\pmod d\}\ =\ 
    \begin{cases}
    0,&n\ \leq\ r,\\[1mm]
    1+\left\lfloor\dfrac{n-1-r}{d}\right\rfloor,&n\ \geq\ r+1.
    \end{cases}
\end{equation}

Let $a=(a_0,\ldots,a_{d-1})\in\Z^d$ be the target point and fix an active coordinate $r\in I(a)$.  Write $\delta_r=\delta_r(a)>0$.  Let $\xi_{r,j}=1$ when the $j\textsuperscript{th}$ update of coordinate $r$ is directed toward the target and $\xi_{r,j}=0$ otherwise.  Depending on which side of the initial cube contains the target, $\xi_{r,j}$ is either $c_{r+1+jd}$ or $1-c_{r+1+jd}$ from the global placement sequence.  Thus the coordinate variables are independent Bernoulli variables with parameter $1/2$.  After $n$ completed moves, define
\begin{equation}
\label{eq:random-coordinate-subset}
E_{r,n}\ :=\ \{w_{r,j}:0\ \leq\ j\ <\ u_r(n),\ \xi_{r,j}\ =\ 1\}
\subseteq\mathcal W_{r,u_r(n)}.
\end{equation}
Then the target has been covered in coordinate $r$ after $n$ moves exactly when
\begin{equation}
\label{eq:coordinate-covering-event}
\sum_{w\in E_{r,n}}w\ \geq\ \delta_r.
\end{equation}

For $E\subseteq\mathcal W_{r,m}$, define
\begin{equation}
\label{eq:subset-membership-word}
x(E)\ =\ (x_{m-1},\ldots,x_0)\in\{0,1\}^m,
\qquad
x_j\ =\ \begin{cases}
1,&w_{r,j}\in E,\\
0,&w_{r,j}\notin E.
\end{cases}
\end{equation}
Its weighted value is
\begin{equation}
\label{eq:weighted-word-value}
V_r(x)\ :=\ \sum_{j=0}^{m-1}x_jw_{r,j}\ =\ \sum_{w\in E}w.
\end{equation}

\subsection{Binary order and the canonical threshold}

For \(x=(x_{m-1},\ldots,x_0)\in\{0,1\}^m\), define
\begin{equation}
\label{eq:binary-rank-definition}
\operatorname{rank}_2(x)\ :=\ \sum_{j=0}^{m-1}x_j2^j.
\end{equation}
For distinct words \(x,y\), let
\(k(x,y)=\max\{j:x_j\neq y_j\}\). We write \(x<_{\mathrm{MSB}}y\)
when \(x_{k(x,y)}<y_{k(x,y)}\).

The superincreasing property, standard in subset-sum problems \cite[Definition~8.34 and Algorithm~8.35]{MVV}, means that the largest weight on which two subsets differ outweighs every possible contribution from smaller weights. Thus this single digit determines the order of their weighted sums.

\begin{theorem}[Binary ordering and rank]
\leavevmode\label{thm:msb-order}
For every \(m\geq1\) and \(x,y\in\{0,1\}^m\),
\begin{equation}
\label{eq:binary-order-equivalence}
V_r(x)\ <\ V_r(y)
\quad\Longleftrightarrow\quad
x\ <_{\mathrm{MSB}\ }y
\quad\Longleftrightarrow\quad
\operatorname{rank}_2(x)\ <\operatorname \ {rank}_2(y).
\end{equation}
Consequently,
\begin{equation}
\label{eq:binary-predecessor-count}
\#\{y\in\{0,1\}^m:V_r(y)\ <\ V_r(x)\}\ =\ \operatorname{rank}_2(x).
\end{equation}
Replacing \(<\) by \(\leq\) increases this count by one.
\end{theorem}

\begin{proof}
For distinct \(x,y\), let \(k=k(x,y)\). If \(x_k=0\) and \(y_k=1\), superincreasingness gives
\begin{equation}
\label{eq:binary-order-proof}
V_r(y)-V_r(x)\ \geq\ w_{r,k}-W_{r,k}\ >\ 0.
\end{equation}
Here the case \(k=0\) uses \(W_{r,0}=0<w_{r,0}\). Interchanging \(x,y\) handles the other case. The same largest differing digit determines the order of the binary ranks. Since these ranks range bijectively over \(0,\ldots,2^m-1\), exactly
\(\operatorname{rank}_2(x)\) words precede \(x\).
\end{proof}

For the positive coordinate distance \(\delta_r\), define
\begin{equation}
\label{eq:coordinate-threshold}
m_r\ :=\ \min\{m\ \geq\ 1:W_{r,m}\ \geq\ \delta_r\}.
\end{equation}
This is the first coordinate-update count at which coverage is possible.

\begin{definition}
\leavevmode\label{def:canonical-threshold}
The canonical threshold word and its ranks are
\begin{equation}
\label{eq:canonical-threshold-ranks}
x_r^\star\ :=\ \min_{<_{\mathrm{MSB}}}
\{x\in\{0,1\}^{m_r}:V_r(x)\ \geq\ \delta_r\},
\qquad
R_r\ :=\ \operatorname{rank}_2(x_r^\star),
\qquad
\kappa_r\ :=\ R_r2^{-m_r}.
\end{equation}
Write \(x_r^\star=(x_{r,m_r-1}^\star,\ldots,x_{r,0}^\star)\).
\end{definition}

Minimality of \(m_r\) gives \(W_{r,m_r-1}<\delta_r\), so the leading digit of \(x_r^\star\) is \(1\). Therefore
\begin{equation}
\label{eq:threshold-rank-bounds}
2^{m_r-1}\ \leq\ R_r\ <\ 2^{m_r},
\qquad
\frac12\ \leq\ \kappa_r\ <\ 1.
\end{equation}

\begin{remark}[Greedy construction]
\leavevmode\label{rem:greedy-threshold}
    Choose the digits from largest to smallest. After fixing the digits above
    \(k\), let
    \begin{equation}
\label{eq:greedy-prefix-weight}
v_{r,k}\ :=\ \sum_{j=k+1}^{m_r-1}x_{r,j}^\star w_{r,j}.
\end{equation}
    Then choose
    \begin{equation}
\label{eq:greedy-digit-rule}
x_{r,k}^\star\ =\ 
    \begin{cases}
    0,&v_{r,k}+W_{r,k}\ \geq\ \delta_r,\\
    1,&v_{r,k}+W_{r,k}\ <\ \delta_r.
    \end{cases}
\end{equation}
    At each stage, the chosen prefix has at least one successful completion. This holds initially because \(W_{r,m_r}\geq\delta_r\). A zero digit is chosen exactly when such a completion remains possible; otherwise every successful completion must use a one. Thus feasibility is preserved, and choosing zero whenever possible produces the least successful word in MSB order.
\end{remark}

\subsection{The exact coordinate covering law}

Fix \(r\in I(a)\), and define
\begin{equation}
\label{eq:coordinate-time-definitions}
\tau_r\ :=\ \inf\left\{m\ \geq\ 0:
\sum_{j=0}^{m-1}\xi_{r,j}w_{r,j}\ \geq\ \delta_r\right\},
\qquad
b_r\ :=\ r+1+d(m_r-1).
\end{equation}
Here \(\tau_r\) counts updates of coordinate \(r\), whereas \(T_r\)
denotes the first global move at which that coordinate is covered.

The key observation is that every weight appearing after the first
feasible update count already exceeds the target distance. Continued
failure therefore requires every newly exposed bit to be zero.

\begin{theorem}[Exact coordinate covering law]
\leavevmode\label{thm:coordinate-covering-law}
For every integer \(m\geq0\),
\begin{equation}
\label{eq:failing-subset-count}
\#\left\{E\subseteq\mathcal W_{r,m}:
\sum_{w\in E}w\ <\ \delta_r\right\}\ =\ 
\begin{cases}
2^m,&m\ <\ m_r,\\
R_r,&m\ \geq\ m_r.
\end{cases}
\end{equation}
For independent fair placement coins, \(\tau_r\) is finite almost surely.
Writing \(G_r:=\tau_r-m_r\), one has
\begin{equation}
\label{eq:coordinate-global-time}
T_r\ =\ r+1+d(\tau_r-1)\ =\ b_r+dG_r
\end{equation}
and
\begin{equation}
\label{eq:coordinate-delay-tail}
\mathbb P(G_r\ >\ k)\ =\ \kappa_r2^{-k}
\qquad(k\ =\ 0,1,\ldots).
\end{equation}
Equivalently,
\begin{equation}
\label{eq:coordinate-delay-mass}
\mathbb P(G_r\ =\ 0)\ =\ 1-\kappa_r,
\qquad
\mathbb P(G_r\ =\ j)\ =\ \kappa_r2^{-j}
\quad(j\ \geq\ 1).
\end{equation}
\end{theorem}

\begin{proof}
If \(m<m_r\), even the sum of all available weights is less than
\(\delta_r\), so every subset fails. At \(m=m_r\),
Theorem~\ref{thm:msb-order} shows that exactly \(R_r\) words precede
the first successful word \(x_r^\star\).

If \(m>m_r\), every new weight satisfies
\begin{equation}
\label{eq:later-weight-target-bound}
w_{r,j}\ >\ W_{r,j}\ \geq\ W_{r,m_r}\ \geq\ \delta_r
\qquad(j\ \geq\ m_r).
\end{equation}
A failing word must therefore have zero in every new position. Its first
\(m_r\) digits must be one of the same \(R_r\) failing words, so the
failing count stays equal to \(R_r\).

All \(2^m\) words are equally likely. Hence
\begin{equation}
\label{eq:coordinate-tail-count}
\mathbb P(\tau_r\ >\ m_r+k)\ =\ R_r2^{-(m_r+k)}\ =\ \kappa_r2^{-k}.
\end{equation}
This tends to zero, proving almost-sure finiteness. It also proves
\eqref{eq:coordinate-delay-tail}; successive tail differences give the
displayed mass function. Finally, the \(m\)\textsuperscript{th} update of coordinate \(r\)
occurs at move \(r+1+d(m-1)\), proving
\eqref{eq:coordinate-global-time}.
\end{proof}

In particular, if \(u_r(n)\) is the number of coordinate-\(r\) updates
among the first \(n\) moves, then
\begin{equation}
\label{eq:coordinate-survival}
\mathbb P(T_r\ >\ n)\ =\ 
\begin{cases}
1,&u_r(n)\ <\ m_r,\\
R_r2^{-u_r(n)},&u_r(n)\ \geq\ m_r.
\end{cases}
\end{equation}

\begin{corollary}[Coordinate moments]
\leavevmode\label{cor:coordinate-moments}
The coordinate-update and global-move moments are
\begin{equation}
\label{eq:coordinate-moment-formulas}
\begin{aligned}
\mathbb E[\tau_r]&\ =\ m_r+2\kappa_r,
&
\operatorname{Var}(\tau_r)&\ =\ 6\kappa_r-4\kappa_r^2,\\
\mathbb E[T_r]&\ =\ b_r+2d\kappa_r,
&
\operatorname{Var}(T_r)&\ =\ d^2(6\kappa_r-4\kappa_r^2).
\end{aligned}
\end{equation}
Moreover, \(2\leq\operatorname{Var}(\tau_r)\leq9/4\).
\end{corollary}

\begin{proof}
The tail in \eqref{eq:coordinate-delay-tail} gives
\(\mathbb E[G_r]=2\kappa_r\) and
\(\mathbb E[G_r^2]=6\kappa_r\). Apply
\(\tau_r=m_r+G_r\) and \(T_r=b_r+dG_r\).
The variance bounds follow by considering
\(6\kappa-4\kappa^2\) on \(1/2\leq\kappa<1\).
\end{proof}

\section{Covering Time in \texorpdfstring{$d$}{d} Dimensions}
\label{sec:covering-times}
We now combine the coordinate subset-sum formulas to obtain the covering-time distribution for a target point
\begin{equation}
\label{eq:target-vector}
a\ =\ (a_0,a_1,\ldots,a_{d-1})\in\Z^d.
\end{equation}

Recall that
\begin{equation}
\label{eq:active-set-recall}
\delta_r\ =\ \delta_r(a)\ =\ \max\{a_r-1,-a_r,0\},
\qquad
I\ =\ I(a)\ =\ \{r:\delta_r\ >\ 0\}.
\end{equation}
Coordinates outside $I$ are already covered by the initial hypercube; for them set $T_r=0$.  We assume below that $I\neq\varnothing$; if $I=\varnothing$, then $T(a)=0$.

For $r\in I$, the update count, feasibility threshold, and coordinate survival law are given by \eqref{eq:update-count}, \eqref{eq:coordinate-threshold}, and \eqref{eq:coordinate-survival}, respectively.

\subsection{Coordinate independence and the exact distribution}

Let $c_k$ be the placement bit used at move $k\geq1$. Coordinate $r$ depends only on the flips
\begin{equation}
\label{eq:coordinate-coin-subsequence}
c_{r+1},\ c_{r+1+d},\ c_{r+1+2d},\ \ldots.
\end{equation}

\begin{lemma}
\leavevmode\label{lemma:coordinate-independence}
The random variables representing the covering times,
\begin{equation}
\label{eq:independent-coordinate-times}
T_0,\ T_1,\ \ldots\ ,\ T_{d-1}
\end{equation}
are mutually independent.
\end{lemma}

\begin{proof}
For each $r$, the random variable $T_r$ is measurable with respect to the sigma-algebra
\begin{equation}
\label{eq:coordinate-sigma-algebra}
\mathcal{G}_r\ =\ \sigma(c_{r+1+kd}:k\ \geq\ 0).
\end{equation}
The index sets
\begin{equation}
\label{eq:coordinate-disjoint-indices}
\{r+1+kd:k\ \geq\ 0\},
\qquad 0\ \leq\ r\ <\ d,
\end{equation}
are pairwise disjoint.  Since the coin flips are mutually independent, the sigma-algebras $\mathcal{G}_0,\ldots,\mathcal{G}_{d-1}$ are mutually independent.  Hence the covering times are mutually independent.
\end{proof}

\begin{definition}
\leavevmode\label{def:full-covering-time}
The target point $a$ is covered when it has been covered in every coordinate. Its covering time is therefore
\begin{equation}
\label{eq:full-covering-maximum}
T\ =\ T(a)\ =\ \max_{0\leq r<d}T_r.
\end{equation}
\end{definition}

\begin{theorem}
\leavevmode\label{thm:exact-distribution}
For every $n\geq0$,
\begin{equation}
\label{eq:full-covering-cdf}
\mathbb{P}(T\ \leq\ n)\ =\ \prod_{r\in I}\bigl(1-\overline F_r(n)\bigr).
\end{equation}
Equivalently,
\begin{equation}
\label{eq:full-covering-survival}
{
\mathbb{P}(T\ >\ n)\ =\ 1-\prod_{r\in I}\bigl(1-\overline F_r(n)\bigr),
}
\end{equation}
where
\begin{equation}
\label{eq:coordinate-survival-function}
\overline F_r(n)\ =\ 
\begin{cases}
1,&u_r(n)\ <\ m_r,\\[1mm]
\dfrac{R_r}{2^{u_r(n)}},&u_r(n)\ \geq\ m_r.
\end{cases}
\end{equation}
\end{theorem}

\begin{proof}
The event $\{T\leq n\}$ is the intersection
\begin{equation}
\label{eq:covering-event-intersection}
\bigcap_{r\in I}\{T_r\ \leq\ n\}.
\end{equation}
Lemma~\ref{lemma:coordinate-independence} therefore gives
\begin{equation}
\label{eq:full-cdf-proof}
\mathbb{P}(T\ \leq\ n)\ =\ \prod_{r\in I}\mathbb{P}(T_r\ \leq\ n)\ =\ \prod_{r\in I}\bigl(1-\overline F_r(n)\bigr).
\end{equation}

Taking complements proves the survival formula.
\end{proof}

\subsection{Cycle decomposition of the geometric tail}

Let
\begin{equation}
\label{eq:maximal-threshold}
m_*\ :=\ \max_{r\in I}m_r.
\end{equation}
Write the number of completed moves as $n=dk+s$, where $0\leq s<d$. Then
\begin{equation}
\label{eq:cycle-update-count}
u_r(dk+s)\ =\ k+\mathbf{1}_{\{r<s\}}.
\end{equation}
In particular, if $k\geq m_*$, every active coordinate is in its geometric regime and
\begin{equation}
\label{eq:cycle-coordinate-tail}
\overline F_r(dk+s)\ =\ \frac{R_r}{2^{k+\mathbf{1}_{\{r<s\}}}}.
\end{equation}

For a nonempty subset $J\subseteq I$, define
\begin{equation}
\label{eq:subset-phase-data}
R_J\ :=\ \prod_{r\in J}R_r,
\qquad
\eta_s(J)\ :=\ \#\{r\in J:r\ <\ s\}.
\end{equation}

\begin{proposition}
\leavevmode\label{prop:inclusionexclusion}
For $k\geq m_*$ and $0\leq s<d$,
\begin{equation}
\label{eq:cycle-inclusion-exclusion}
\mathbb{P}(T\ >\ dk+s)\ =\ 
\sum_{\varnothing\neq J\subseteq I}
(-1)^{|J|+1}
R_J2^{-k|J|-\eta_s(J)}.
\end{equation}

\end{proposition}

\begin{proof}
Substitute the geometric coordinate tails into
\begin{equation}
\label{eq:product-survival-expansion}
1-\prod_{r\in I}(1-\overline F_r(dk+s))
\end{equation}
and expand the product by inclusion-exclusion.  For a fixed $J$, the product
of its coordinate tails is
\begin{equation}
\label{eq:subset-tail-product}
\prod_{r\in J}
\frac{R_r}{2^{k+\mathbf{1}_{\{r<s\}}}}\ =\ R_J2^{-k|J|-\eta_s(J)}.
\end{equation}
\end{proof}

To shorten the moment formulas, define the phase sums
\begin{equation}
\label{eq:phase-sums}
H_J\ =\ \sum_{s=0}^{d-1}2^{-\eta_s(J)}
\quad\text{ and }\quad
K_J\ =\ \sum_{s=0}^{d-1}(2s+1)2^{-\eta_s(J)}.
\end{equation}
Both depend only on the positions of the coordinates in $J$ within the cyclic
update order.

\subsection{Moment Formulas}

For a nonnegative integer-valued random variable,
\begin{equation}
\label{eq:expectation-tail-sum}
\mathbb{E}[T]\ =\ \sum_{n=0}^{\infty}\mathbb{P}(T\ >\ n).
\end{equation}
We split this sum immediately before the first complete geometric cycle as
\begin{equation}
\label{eq:expectation-split}
\mathbb{E}[T]\ =\ E_{\mathrm{fin}}+E_{\mathrm{tail}},
\end{equation}
where
\begin{equation}
\label{eq:expectation-finite-part}
E_{\mathrm{fin}}\ =\ \sum_{n=0}^{dm_*-1}\mathbb{P}(T\ >\ n).
\end{equation}

\begin{theorem}[Exact expected covering time]
\leavevmode\label{thm:exact-expected}
Fix $d\geq2$ and a lattice target $a$ with $I(a)\neq\varnothing$, and assume independent fair placement coins. With the quantities defined above, the expected covering time is
\begin{equation}
\label{eq:expected-covering-time}
{
\begin{aligned}
\mathbb{E}[T]\ =\ {}&
\sum_{n=0}^{dm_*-1}
\left[
1-\prod_{r\in I}\bigl(1-\overline F_r(n)\bigr)
\right]\\
&+
\sum_{\varnothing\neq J\subseteq I}
(-1)^{|J|+1}
\frac{R_JH_J\,2^{-m_*|J|}}
{1-2^{-|J|}}.
\end{aligned}
}
\end{equation}
\end{theorem}

\begin{proof}
Only the tail requires evaluation.  Proposition~\ref{prop:inclusionexclusion} gives
\begin{equation}
\label{eq:expectation-tail-evaluation}
\begin{aligned}
E_{\mathrm{tail}}
&\ =\ \sum_{k=m_*}^{\infty}\sum_{s=0}^{d-1}
\mathbb{P}(T\ >\ dk+s)\\
&\ =\ 
\sum_{\varnothing\neq J\subseteq I}
(-1)^{|J|+1}R_J
\sum_{s=0}^{d-1}2^{-\eta_s(J)}
\sum_{k=m_*}^{\infty}2^{-k|J|}.
\end{aligned}
\end{equation}
Since
\begin{equation}
\label{eq:geometric-tail-series}
\sum_{k=m_*}^{\infty}2^{-k|J|}\ =\ \frac{2^{-m_*|J|}}{1-2^{-|J|}},
\end{equation}
the result follows.
\end{proof}

Although the original tail is an infinite series, the theorem expresses it as
a sum over the $2^{|I|}-1$ nonempty active-coordinate subsets.  Thus, for fixed $d$, the
only unsummed part is the finite transient range $0\leq n<dm_*$.

The second tail-sum identity is
\begin{equation}
\label{eq:second-moment-tail-sum}
\mathbb{E}[T^2]\ =\ \sum_{n=0}^{\infty}(2n+1)\mathbb{P}(T\ >\ n).
\end{equation}
Define
\begin{equation}
\label{eq:second-moment-finite-part}
M_{2,\mathrm{fin}}\ =\ \sum_{n=0}^{dm_*-1}(2n+1)
\left[
1-\prod_{r\in I}\bigl(1-\overline F_r(n)\bigr)
\right].
\end{equation}

\begin{theorem}[Exact second moment and variance]
\leavevmode\label{thm:exact-second-moment}
Fix $d\geq2$ and a lattice target $a$ with $I(a)\neq\varnothing$, and assume independent fair placement coins. For each nonempty $J\subseteq I$, let
\begin{equation}
\label{eq:subset-geometric-ratio}
z_J\ =\ 2^{-|J|}.
\end{equation}
Then
\begin{equation}
\label{eq:second-moment-formula}
{
\begin{aligned}
\mathbb{E}[T^2]\ =\ {}&M_{2,\mathrm{fin}}\\
&+
\sum_{\varnothing\neq J\subseteq I}
(-1)^{|J|+1}R_Jz_J^{m_*}
\left[
\frac{2dm_* H_J+K_J}{1-z_J}
+\frac{2dz_JH_J}{(1-z_J)^2}
\right].
\end{aligned}
}
\end{equation}
Consequently,
\begin{equation}
\label{eq:variance-formula}
{
\operatorname{Var}(T)\ =\ \mathbb{E}[T^2]-\mathbb{E}[T]^2,
}
\end{equation}
where $\mathbb{E}[T]$ is given by Theorem~\ref{thm:exact-expected}.
\end{theorem}

\begin{proof}
For a fixed nonempty set $J$, its contribution to the second-moment tail is
\begin{equation}
\label{eq:second-moment-subset-contribution}
\begin{aligned}
&R_J\sum_{s=0}^{d-1}2^{-\eta_s(J)}
\sum_{k=m_*}^{\infty}(2dk+2s+1)z_J^k.
\end{aligned}
\end{equation}
The geometric identities
\begin{equation}
\label{eq:geometric-series}
\sum_{k=m_*}^{\infty}z^k\ =\ \frac{z^{m_*}}{1-z}
\end{equation}
and
\begin{equation}
\label{eq:shifted-geometric-first-moment}
\sum_{k=m_*}^{\infty}kz^k\ =\ z^{m_*}\left(
\frac{m_*}{1-z}+\frac{z}{(1-z)^2}
\right)
\end{equation}
give
\begin{equation}
\label{eq:second-moment-subset-evaluation}
R_Jz_J^{m_*}
\left[
\frac{2dm_*H_J+K_J}{1-z_J}
+\frac{2dz_JH_J}{(1-z_J)^2}
\right].
\end{equation}
Summing with the inclusion--exclusion signs and adding the finite transient
part proves the result.
\end{proof}

\section{Limiting Behavior}
\label{sec:limiting-behavior}
Theorems~\ref{thm:exact-distribution}, \ref{thm:exact-expected}, and \ref{thm:exact-second-moment} of Section~\ref{sec:covering-times} determine the covering-time distribution for every fixed target.  We now use them to describe the behavior of the covering time as the target moves away from the initial hypercube.  We show that the expected covering time grows logarithmically, the excess beyond its deterministic threshold has a uniformly exponential tail, and all centered moments remain bounded when the dimension is fixed.

Throughout this section, $d\geq2$ is fixed and $I(a)\neq\varnothing$. For $a\in\Z^d$, recall
\begin{equation}
\label{eq:limiting-target-distances}
\delta_r(a)\ =\ \max\{a_r-1,-a_r,0\},
\qquad
I(a)\ =\ \{r:\delta_r(a)\ >\ 0\},
\end{equation}
and write
\begin{equation}
\label{eq:limiting-distance}
\delta_\infty(a)\ :=\ \max_{0\leq r<d}\delta_r(a).
\end{equation}
When the target is fixed, we abbreviate $\delta_r=\delta_r(a)$ and $I=I(a)$.
For $r\in I(a)$, let $b_r=r+1+d(m_r-1)$ be the earliest move at which coordinate $r$ can be covered, and let
\begin{equation}
\label{eq:earliest-full-covering}
B(a)\ :=\ \max_{r\in I(a)}b_r
\end{equation}
be the earliest move at which all active coordinates can be covered.

\subsection{Deterministic location}

For each active coordinate, minimality of \(m_r\) gives
\begin{equation}
\label{eq:threshold-bracketing}
W_{r,m_r-1}\ <\ \delta_r\ \leq\ W_{r,m_r}.
\end{equation}

\begin{corollary}
\leavevmode\label{cor:threshold-asymptotics}
For fixed \(d\geq2\), uniformly over active coordinates and lattice targets,
\begin{equation}
\label{eq:threshold-asymptotics}
dm_r\ =\ \log_{\rho_d}\delta_r+O_d(1),
\qquad
B(a)\ =\ \log_{\rho_d}\delta_\infty(a)+O_d(1).
\end{equation}
\end{corollary}

\begin{proof}
Proposition~\ref{prop:coordinate-superincreasing} gives
\begin{equation}
\label{eq:cumulative-weight-bound}
w_{r,m-1}\ \leq\ W_{r,m}\ <\ 2w_{r,m-1}
\qquad(m\ \geq\ 1).
\end{equation}
Since \(w_{r,m-1}=F_{r+dm-1}^{(d)}\), Theorem~\ref{theorem:asymptotic} implies
\begin{equation}
\label{eq:cumulative-weight-growth}
\log_{\rho_d}W_{r,m}\ =\ r+dm-1+O_d(1).
\end{equation}
For \(m_r\geq2\), apply this to the two partial sums bracketing
\(\delta_r\). If \(m_r=1\), then \(1\leq\delta_r\leq w_{r,0}\),
so the same estimate holds after enlarging the constant. Thus \(b_r=r+1+d(m_r-1)=\log_{\rho_d}\delta_r+O_d(1)\). Taking the maximum over the finitely many active coordinates proves the assertion for \(B(a)\).
\end{proof}

\subsection{Uniform fluctuations and their consequences}

Assume \(I(a)\neq\varnothing\). By Theorem~\ref{thm:coordinate-covering-law},
\(T_r=b_r+dG_r\) and
\(\mathbb P(G_r>k)=\kappa_r2^{-k}\).
Set \(G_{\max}:=\max_{r\in I(a)}G_r\).
Since \(T=\max_{r\in I(a)}T_r\) and \(B(a)=\max_{r\in I(a)}b_r\),
\begin{equation}
\label{eq:limiting-delay-comparison}
0\ \leq\ T-B(a)\ =\ \max_{r\in I(a)}\{b_r-B(a)+dG_r\}\ \leq\ dG_{\max}.
\end{equation}
Thus the fluctuations are controlled by at most \(d\) coordinate delays, each having a geometric tail independent of the scale of the target.

\begin{theorem}
\leavevmode\label{thm:uniform-fluctuations}
For every target with \(I(a)\neq\varnothing\) and every integer \(k\geq0\),
\begin{equation}
\label{eq:uniform-delay-tail}
\mathbb P(T\ >\ B(a)+dk)\ \leq\ \min\{1,|I(a)|2^{-k}\}\ \leq\ \min\{1,d2^{-k}\}.
\end{equation}
For every fixed real \(q\geq1\),
\begin{equation}
\label{eq:uniform-delay-moments}
\mathbb E[|T-B(a)|^q]\ =\ O_{d,q}(1),
\end{equation}
uniformly over the target.
\end{theorem}

\begin{proof}
The union bound and \eqref{eq:threshold-rank-bounds} give
\begin{equation}
\label{eq:maximum-delay-union-bound}
\mathbb P(G_{\max}\ >\ k)\ \leq\ \sum_{r\in I(a)}\mathbb P(G_r\ >\ k)\ =\ \sum_{r\in I(a)}\kappa_r2^{-k}\ \leq\ |I(a)|2^{-k}.
\end{equation}
Together with \eqref{eq:limiting-delay-comparison}, this proves the concentration bound. For the moments, the tail-sum identity gives
\begin{equation}
\label{eq:maximum-delay-moments}
\begin{aligned}
\mathbb E[G_{\max}^q]
&\ =\ \sum_{k=0}^{\infty}\bigl((k+1)^q-k^q\bigr)
   \mathbb P(G_{\max}\ >\ k)\\
&\ \leq\ \sum_{k=0}^{\infty}\bigl((k+1)^q-k^q\bigr)
   \min\{1,d2^{-k}\}\ =\ O_{d,q}(1).
\end{aligned}
\end{equation}
The series converges because its summands are bounded by an exponential times a polynomial. Apply \eqref{eq:limiting-delay-comparison} once more.
\end{proof}

\begin{corollary}
\leavevmode\label{cor:covering-asymptotics}
For fixed \(d\),
\begin{equation}
\label{eq:covering-asymptotics}
\mathbb E[T]\ =\ B(a)+O_d(1)\ =\ \log_{\rho_d}\delta_\infty(a)+O_d(1),
\qquad
\operatorname{Var}(T)\ =\ O_d(1),
\end{equation}
uniformly over targets with \(I(a)\neq\varnothing\).
As \(\delta_\infty(a)\to\infty\),
\begin{equation}
\label{eq:relative-lq-convergence}
\frac{T}{\log_{\rho_d}\delta_\infty(a)}\longrightarrow1
\quad\text{in }L^q
\end{equation}
for every fixed \(q\geq1\). Moreover,
\begin{equation}
\label{eq:coefficient-of-variation}
\frac{\sqrt{\operatorname{Var}(T)}}{\mathbb E[T]}
\longrightarrow0.
\end{equation}
\end{corollary}

\begin{proof}
The \(q=1\) case of Theorem~\ref{thm:uniform-fluctuations} gives
\(\mathbb E[T]-B(a)=O_d(1)\), and the \(q=2\) case gives
\begin{equation}
\label{eq:uniform-variance-bound}
\operatorname{Var}(T)\ \leq\ \mathbb E[(T-B(a))^2]\ =\ O_d(1).
\end{equation}
Use Corollary~\ref{cor:threshold-asymptotics} for the logarithmic expression. For any fixed \(q\geq1\), the same two results imply
\begin{equation}
\label{eq:logarithmic-centered-moments}
\begin{aligned}
\mathbb E\bigl[|T-\log_{\rho_d}\delta_\infty(a)|^q\bigr]\ \leq\ 2^{q-1}\bigl(
&\mathbb E[|T-B(a)|^q]\\
&+|B(a)-\log_{\rho_d}\delta_\infty(a)|^q
\bigr)\ =\ O_{d,q}(1).
\end{aligned}
\end{equation}
Dividing by \((\log_{\rho_d}\delta_\infty(a))^q\) proves \(L^q\) convergence. The final conclusion follows because the variance is bounded and the expectation tends to infinity.
\end{proof}

\subsection{The centered covering-time law}

Theorem~\ref{thm:uniform-fluctuations} implies that
\begin{equation}
\label{eq:tight-centered-family}
\{T(a)-B(a):a\in\mathbb Z^d,\ I(a)\ \neq\ \varnothing\}
\end{equation}
is tight for fixed \(d\): its members are nonnegative and their tails are uniformly bounded by \(d2^{-k}\) at \(dk\).

The independent coordinate delays also give the exact identity
\begin{equation}
\label{eq:maximum-delay-exact-tail}
\mathbb P(G_{\max}\ >\ k)\ =\ 1-\prod_{r\in I(a)}(1-\kappa_r2^{-k}).
\end{equation}
Indeed, the coordinate delays depend on disjoint sets of independent placement coins. This identity is a law for \(G_{\max}\), not generally for \((T-B(a))/d\), because
\begin{equation}
\label{eq:centered-covering-identity}
T-B(a)\ =\ \max_{r\in I(a)}\{b_r-B(a)+dG_r\}.
\end{equation}
Consequently, the centered law depends on the active-coordinate set, the ranks \(\kappa_r\), and the offsets \(b_r-B(a)\). Tightness guarantees weakly convergent subsequences, but does not by itself identify a universal centered limiting law.

\subsection{The planar case}

When $d=2$, coordinate $0$ uses the odd-indexed Fibonacci numbers $F_{1+2j}$, coordinate $1$ uses the positive even-indexed Fibonacci numbers $F_{2+2j}$, and
\begin{equation}
\label{eq:planar-dominant-root}
\rho_2\ =\ \varphi\ =\ \frac{1+\sqrt5}{2}.
\end{equation}

\begin{corollary}
\leavevmode\label{cor:planar-covering}
For the planar random Fibonacci tiling and a target with $I(a)\neq\varnothing$,
\begin{equation}
\label{eq:planar-covering-asymptotics}
\mathbb{E}[T]\ =\ \log_\varphi\delta_\infty(a)+O(1),
\qquad
\operatorname{Var}(T)\ =\ O(1).
\end{equation}
If $B(a)$ denotes the earliest move at which both coordinates can have been covered, then
\begin{equation}
\label{eq:planar-delay-tail}
\mathbb{P}\bigl(T\ >\ B(a)+2k\bigr)\ \leq\ \min\{1,2^{1-k}\}
\qquad(k\ \geq\ 0).
\end{equation}
Moreover,
\begin{equation}
\label{eq:planar-relative-convergence}
\frac{T}{\log_\varphi\delta_\infty(a)}
\longrightarrow1
\end{equation}
in $L^q$ for every fixed $q\geq1$ as $\delta_\infty(a)\to\infty$.
\end{corollary}

\begin{remark}[Parity effects]
\leavevmode\label{rem:parity-effects}
The exact expectation and variance formulas retain the distinction between the coordinate updated on even global steps and the coordinate updated on odd global steps through their phase factors.  These parity corrections are bounded, so they are absorbed into the $O(1)$ terms in the limiting formulas.
\end{remark}

\subsection{Dependence on the dimension}

We finally allow $d\to\infty$ to examine the coefficient of $\log \delta_\infty(a)$ in the deterministic scale. Let $\operatorname{W}_0$ denote the principal real branch of the Lambert function \cite[Section~4]{CGHJK}, defined by $\operatorname{W}_0(x)e^{\operatorname{W}_0(x)}=x$ for $x>0$.

\begin{proposition}
\leavevmode\label{prop:limiting-dimension}
As $d\to\infty$,
\begin{equation}
\label{eq:dimension-root-asymptotics}
\rho_d\ =\ 1+\frac{\operatorname{W}_0(d-1)}{d-1}+o(1/d),
\qquad
\log\rho_d\ \sim\ \frac{\operatorname{W}_0(d)}{d}.
\end{equation}
Consequently, the logarithmic scale satisfies
\begin{equation}
\label{eq:dimension-logarithmic-scale}
\log_{\rho_d}\delta_\infty(a)\ \sim\ \frac{d}{\operatorname{W}_0(d)}\log \delta_\infty(a)
\qquad(\delta_\infty(a)\ >\ 1).
\end{equation}
\end{proposition}

\begin{proof}
Write $m=d-1$ and $u=m(\rho_d-1)$. Taking logarithms in
$\rho_d^{d-1}(\rho_d-1)=1$ gives
\begin{equation}
\label{eq:dimension-log-root-equation}
m\log(1+u/m)+\log u\ =\ \log m.
\end{equation}
Since $0<u/m<1$, the inequality $\log(1+x)\geq x/2$ for $0\leq x\leq 1$ implies $u\leq 2\log m$ whenever $u\geq 1$. Thus $u=O(\log m)$, and expanding the logarithm yields
\begin{equation}
\label{eq:dimension-root-error}
u+\log u\ =\ \log m+O\!\left(\frac{(\log m)^2}{m}\right)\ =\ \log m+o(1).
\end{equation}
The function $v\mapsto v+\log v$ has derivative at least $1$ on $(0,\infty)$, and $\operatorname{W}_0(m)+\log \operatorname{W}_0(m)=\log m$. Hence $u=\operatorname{W}_0(m)+o(1)$. Substituting into $\rho_d=1+u/m$ proves the first estimate. The others follow from $\log(1+u/m)\sim u/m$ and $\operatorname{W}_0(d-1)\sim \operatorname{W}_0(d)$.
\end{proof}

This comparison concerns the deterministic growth scale only. The bounds denoted by $O_d$ may depend on $d$, so these results do not by themselves establish a joint limit in which both dimension and target distance grow.


\section{Random Fibonacci Spirals}
\label{sec:random-spirals}

The usual Fibonacci spiral joins quarter-circles drawn in successive
squares of the standard tiling.
In a random tiling, the entrance and exit
corners of a square may be adjacent rather than opposite, so a
quarter-circle no longer suffices. We retain the usual arcs whenever
possible and otherwise use a quintic connector followed by a smaller
quarter-circle. Throughout this section, the resulting curve is
parametrized by arclength, so it is traversed at unit speed.

We fix the splice at half the side length. This choice gives a simple
construction whose segments remain inside their assigned squares and
whose concatenation has no self-intersections. The quintic connector
is based on the construction proposed by L\'aszl\'o Szalay.
The geometric questions at the end of the section discuss other splice
choices, containment, self-intersection, and further extensions.

\begin{figure}[h]
    \centering
    \includegraphics[width=0.47\textwidth]{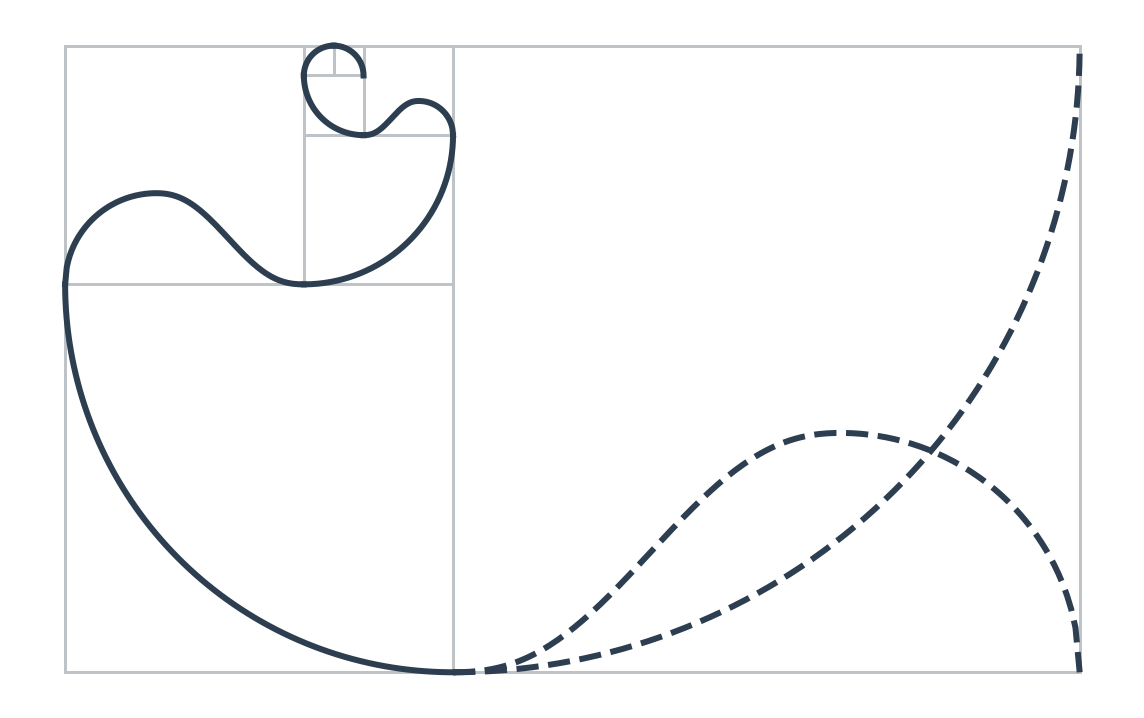}
    \caption{One-step lookahead, shown schematically. Dashed curves are
    alternative continuations determined by the next placement; only
    the selected continuation is included in the spiral.}
    \label{fig:random-spiral}
\end{figure}

\subsection{Entrance and exit corners}

Let $\mathcal Q_0=[0,1]^2$, and let $\mathcal Q_n$ be the square of side length $F_n$
added at move $n\geq1$. Fix a placement word
\begin{equation}
    \omega=(\omega_1,\omega_2,\ldots)\in\{0,1\}^{\mathbb N},
    \label{eq:spiral-placement-word}
\end{equation}
with $\omega_n=c_n$ for a fixed realization of the placement bits
in Section~2.
Write $v_n$ for the unit vector in the direction in which $\mathcal Q_n$ is
attached. Consecutive vectors $v_n$ and $v_{n+1}$ are perpendicular.

If $\mathcal Q_1$ is attached to the right, take $\gamma_0$ to be the
quarter-circle centered at $(1,1)$, traversed from $(0,1)$ to $(1,0)$.
If $\mathcal Q_1$ is attached to the left, reflect this arc across $x=1/2$.
Let $p_0$ be its endpoint. In either case, $p_0$ is a corner of $\mathcal Q_1$
and the terminal tangent of $\gamma_0$ is $v_1$.

For $n\geq1$, let $p_n$ be the unique common corner of $\mathcal Q_n$ and
$\mathcal Q_{n+1}$. We construct a segment $\gamma_n$ from $p_{n-1}$ to $p_n$
with entrance tangent $v_n$ and exit tangent $v_{n+1}$.

\begin{lemma}
\label{lem:spiral-configurations}
For every $n\geq1$, the points $p_{n-1}$ and $p_n$ are corners of
$\mathcal Q_n$, and exactly one of the following holds:
\begin{align}
    p_n-p_{n-1}&=F_n(v_n+v_{n+1}),
        \label{eq:spiral-standard-corners}\\
    p_n-p_{n-1}&=F_nv_n.
        \label{eq:spiral-nonstandard-corners}
\end{align}
In the first case the corners are opposite; in the second they are
adjacent. We call these the \emph{standard} and \emph{nonstandard}
configurations, respectively.
\end{lemma}

\begin{proof}
The square $\mathcal Q_{n+1}$ shares a side of $\mathcal Q_n$ and extends beyond it
toward the older squares. Since $F_{n+1}>F_n$ for $n\geq1$, precisely
one endpoint of the shared side is a corner of both squares. This is
the endpoint on the outer face newly created by attaching $\mathcal Q_n$.
Thus $p_n$ is well-defined and lies on that outer face.

The entrance $p_{n-1}$ lies on the opposite, attaching face of $\mathcal Q_n$.
For $n=1$ this follows from the choice of $p_0$; for $n\geq2$ it
follows from the definition of $p_{n-1}$. The displacement from
$p_{n-1}$ to $p_n$ therefore has component $F_nv_n$. Its perpendicular
component is either zero or $F_nv_{n+1}$, according to which of the
two sides receives $\mathcal Q_{n+1}$.
\end{proof}

In a standard configuration, let $\gamma_n$ be the usual
quarter-circle of radius $F_n$ joining the two corners with the
prescribed tangents. The nonstandard configuration requires the
connector defined next.

\subsection{The local nonstandard segment}

Work first in the reference square $[0,B_n]^2$, where $B_n=F_n$ and
$n\geq1$. The entrance and exit are $(0,0)$ and $(B_n,0)$, with tangent
directions $(1,0)$ and $(0,-1)$. Set
\[
    \mu_n=\frac{F_n}{F_{n-1}},
    \qquad
    \sigma_5(t)=10t^3-15t^4+6t^5,
\]
and define
\begin{equation}
    P_n(x)
    =\frac{B_n}{2}\sigma_5(t)
     +\frac{B_n\mu_n}{8}t^2(1-t)^3
     -\frac{B_n}{4}t^3(1-t)^2,
    \qquad t=\frac{2x}{B_n},\quad 0\leq x\leq\frac{B_n}{2}.
    \label{eq:quintic-polynomial}
\end{equation}
The reference segment follows the graph of $P_n$ and then the
quarter-circle
\begin{equation}
    y=\sqrt{\frac{B_n^2}{4}-\left(x-\frac{B_n}{2}\right)^2},
    \qquad \frac{B_n}{2}\leq x\leq B_n.
    \label{eq:reference-circular-arcs}
\end{equation}
To place this segment in $\mathcal Q_n$, apply the rigid motion
\begin{equation}
    (x,y)\longmapsto p_{n-1}+xv_n-yv_{n+1}.
    \label{eq:spiral-transport}
\end{equation}
By Lemma~\ref{lem:spiral-configurations}, this sends the reference
square onto $\mathcal Q_n$ and gives the required endpoints and tangent
directions.

\begin{theorem}
\label{thm:spiral-quintic}
The polynomial $P_n$ is the unique polynomial of degree at most five
satisfying the endpoint conditions~\eqref{eq:connector-conditions}.
The reference segment is simple, joins $(0,0)$ to $(B_n,0)$ with the
prescribed tangents, and lies strictly inside $[0,B_n]^2$ except at
these endpoints. The quintic and the outgoing quarter-circle meet
with $C^2$ continuity. Moreover, $P_n$ is increasing on $[0,B_n/2]$.
\end{theorem}

\begin{proof}
Direct differentiation gives
\begin{equation}
\begin{aligned}
    P_n(0)&\ =\ 0,& P_n'(0)&\ = \ 0,&P_n''(0)& \ = \ \frac{1}{F_{n-1}},\\
    P_n(B_n/2)& \ = \ B_n/2,&P_n'(B_n/2)& \ = \ 0,&P_n''(B_n/2)& \ =\ -\frac{2}{B_n}.
\end{aligned}
\label{eq:connector-conditions}
\end{equation}
The outgoing circle has the same value and first two derivatives
at $x=B_n/2$, and its terminal tangent is $(0,-1)$. The polynomial
is the quintic Hermite interpolant for these endpoint data
\cite[Proposition~32]{dB}. For uniqueness, the difference of two
solutions is divisible by $x^3(x-B_n/2)^3$ and has degree at most
five, so it vanishes.

For containment and monotonicity, write
\[
    \frac{P_n(Bt/2)}{B_n}
     \ = \ \sum_{j=0}^5 b_j\binom{5}{j}t^j(1-t)^{5-j},
    \qquad
    (b_0,\ldots,b_5)
     \ = \ \left(0,0,\frac{\mu_n}{80},\frac{19}{40},\frac12,\frac12\right).
\]
Since $1\leq\mu_n\leq2$, these coefficients are nondecreasing and
lie in $[0,1/2]$. The Bernstein basis functions are nonnegative
and sum to one. Hence $0<P_n(x)<B_n/2$ for $0<x<B_n/2$.
Differentiating the Bernstein expression gives nonnegative
coefficients $5(b_{j+1}-b_j)$, so $P_n$ is increasing.

Both the connector and the outgoing arc are therefore inside the
square, meeting its boundary only at $(0,0)$ and $(B_n,0)$.
The full reference segment is simple because its $x$ coordinate
is strictly increasing along the traversal.
\end{proof}

\begin{figure}[b]
    \centering
    \includegraphics[width=0.84\textwidth]{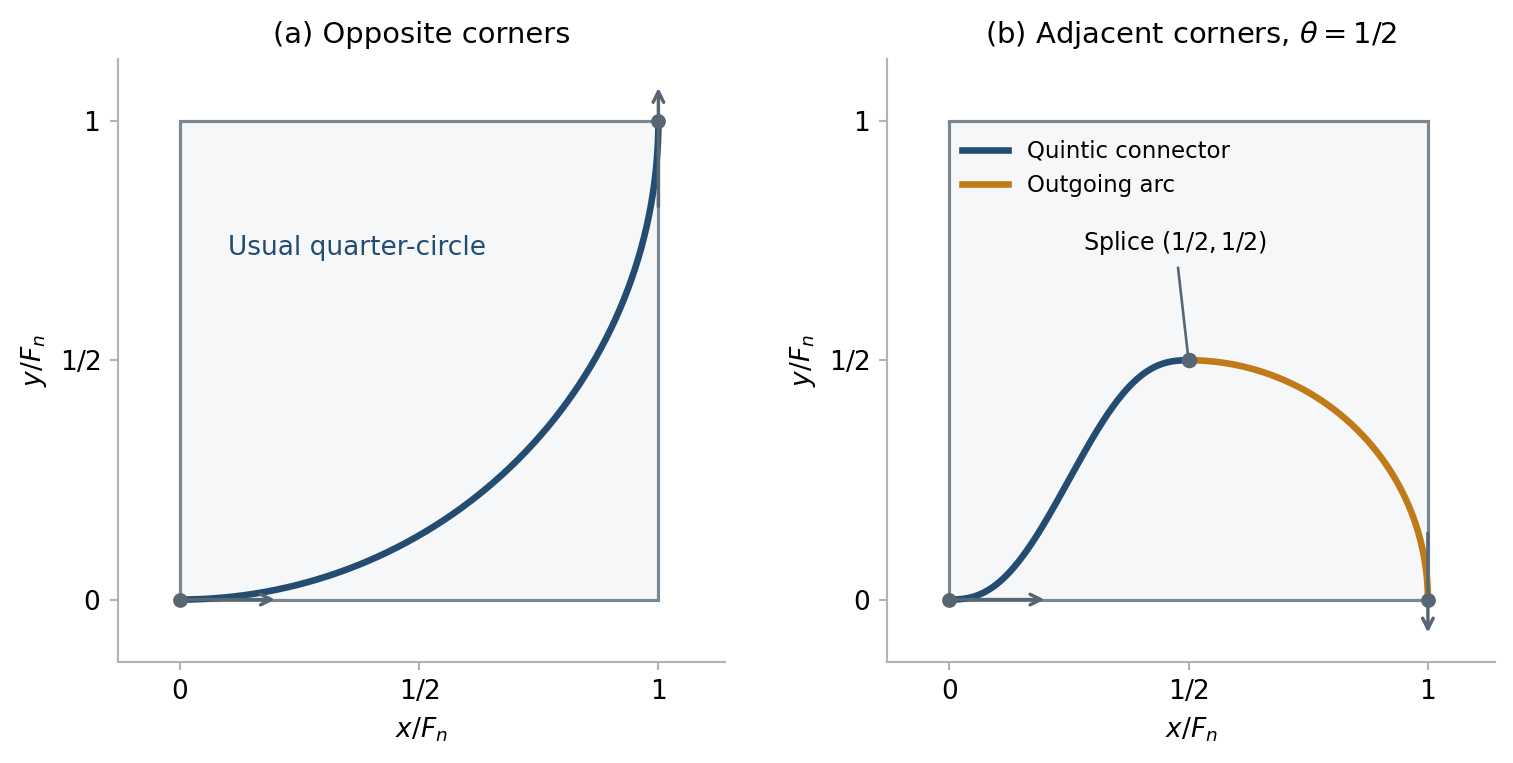}
    \caption{The two local configurations, normalized to a unit square.}
    \label{fig:spiral-local-configurations}
\end{figure}

\subsection{Definition and global regularity}

\begin{definition}
\label{def:random-fib-spiral}
Let $\Gamma(\omega)$ be the curve obtained by following
$\gamma_0,\gamma_1,\gamma_2,\ldots$ in order, using a quarter-circle
in each standard configuration and the transported reference
segment in each nonstandard configuration.
\end{definition}

\begin{remark}
\label{rem:spiral-lookahead}
The segment $\gamma_n$
is determined as soon as $\mathcal Q_{n+1}$ is placed, so the construction
uses one-step lookahead. A word of length $N\geq1$ determines the
curve through $\mathcal Q_{N-1}$. Figure~\ref{fig:random-spiral} illustrates
this dependence schematically; Figure~\ref{fig:spiral-local-configurations}
shows the exact local curves used here.
\end{remark}

\begin{theorem}
\label{thm:spiral-global}
For every placement word $\omega$, the construction defines a
simple $C^1$ curve $\Gamma(\omega)$, with each segment contained
in its assigned square. Within each nonstandard segment, the
quintic and circular arc meet with $C^2$ continuity.

For the standard word $\omega=(0,0,1,1,0,0,1,1,\ldots)$, the
construction recovers the usual Fibonacci spiral.
\end{theorem}

\begin{proof}
Lemma~\ref{lem:spiral-configurations} supplies the two possible
configurations, and Theorem~\ref{thm:spiral-quintic} supplies the
nonstandard segment. At $p_n$, the terminal tangent of $\gamma_n$
and the entrance tangent of $\gamma_{n+1}$ are both $v_{n+1}$.
Their unit-speed parametrizations therefore join with a continuous
derivative. The internal $C^2$ claim follows from
Theorem~\ref{thm:spiral-quintic}.

Each segment lies in the interior of its square apart from its
two endpoints. The earlier segments lie in the rectangle present
before $\mathcal Q_n$ was attached. The interior of $\mathcal Q_n$ lies outside that
rectangle, as does its exit corner $p_n$ on the new outer face.
Consequently, $\gamma_n$ meets the preceding curve only at its
entrance $p_{n-1}$. Since each individual segment is simple,
induction proves that the concatenation has no self-intersections.

Each segment has length at least $F_n$, the minimum distance
between distinct corners of $\mathcal Q_n$. Since $\sum_{n\geq0}F_n=\infty$,
the joining times do not accumulate in finite arclength, and the
curve is defined for all nonnegative arclengths.

For the displayed standard word, the attachment directions repeat
right, up, left, down. Every configuration has opposite entrance
and exit corners, so every segment is the corresponding usual
quarter-circle, beginning with the specified initial arc.
\end{proof}

\begin{remark}
\label{rem:spiral-curvature}
As in the usual Fibonacci spiral, curvature may jump between
consecutive squares; the construction guarantees $C^1$ continuity
at these joins.
\end{remark}

\subsection{Length}

\begin{proposition}
\label{prop:spiral-length-scaling}
If $\ell_n$ is the length of $\gamma_n$, then
\begin{equation}
    F_n \ \leq\ \ ell_n \ \leq\ \left(1+\frac{\pi}{4}\right)F_n
    \qquad(n\geq0).
    \label{eq:spiral-length-bounds}
\end{equation}
Consequently,
\begin{equation}
    F_{N+2}-1
     \ \leq \ \sum_{n=0}^N\ell_n
     \ \leq \ \left(1+\frac{\pi}{4}\right)(F_{N+2}-1),
    \label{eq:cumulative-length-bounds}
\end{equation}
uniformly over the placement word. In particular,
\begin{equation}
    \sum_{n=0}^N\ell_n \ = \ \Theta(F_{N+2}).
    \label{eq:cumulative-spiral-length}
\end{equation}
\end{proposition}

\begin{proof}
The lower bound is the distance between distinct square corners.
The initial arc and every standard segment have length $\pi F_n/2$.
For a nonstandard segment, the connector is increasing in both
coordinates, and its total horizontal and vertical displacements
are each $F_n/2$. Its length is therefore at most $F_n$. Adding
the quarter-circle of radius $F_n/2$ gives the upper bound.
Finally, sum over $n$ and use
$\sum_{n=0}^N F_n=F_{N+2}-1$.
\end{proof}

\subsection{Geometric questions}
\label{section:questions}

The midpoint construction settles containment and self-intersection
for the curve defined above. These remain important questions when
the splice or the matching rule is changed.

For comparison, allow the splice to occur at $x=\theta F_n$, where
$0<\theta<1$, and use an outgoing quarter-circle of radius
$(1-\theta)F_n$. The corresponding normalized quintic is
\begin{equation}
    g_{\theta,\mu}(t)
     \ = \ (1-\theta)\sigma_5(t)
    +\frac{\theta^2\mu}{2}t^2(1-t)^3
    -\frac{\theta^2}{2(1-\theta)}t^3(1-t)^2,
    \qquad 0\leq t\leq1.
    \label{eq:normalized-connector-polynomial}
\end{equation}
The connector is $(x,y)=F_n(\theta t,g_{\theta,\mu_n}(t))$.
Write $\Gamma_\theta(\omega)$ for the ordered curve obtained using
this segment and the same corner rule. At $\theta=1/2$, it reduces
to the construction above.

\begin{figure}[h]
    \centering
    \includegraphics[width=0.87\textwidth]{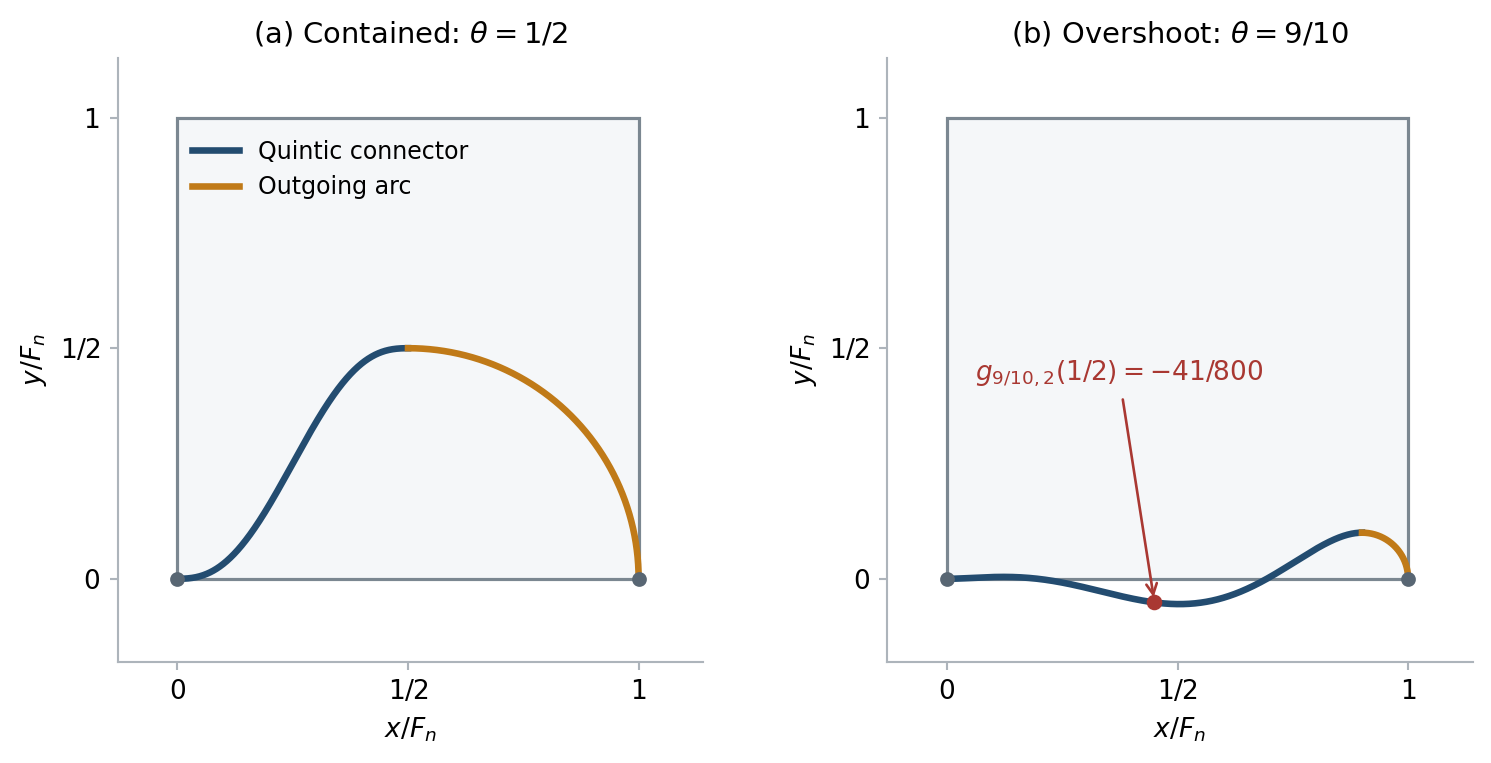}
    \caption{Containment depends on the splice choice. Both plots use
    $\mu=2$. The midpoint segment stays inside its square, whereas
    the connector with $\theta=9/10$ crosses below the lower edge.}
    \label{fig:spiral-containment}
\end{figure}

\begin{enumerate}
\item \textbf{Choice of splice.}
Is the midpoint geometrically distinguished, or is another choice,
such as $\theta=\varphi^{-1}$, more natural? One could instead choose
$\theta$ to minimize length, maximum curvature, or bending energy
\begin{equation}
    \int \mathcal K(s)^2\,ds,
    \label{eq:bending-energy}
\end{equation}
where $\mathcal K$ is curvature and $s$ is arclength.\\

\item \textbf{Containment.}
For which values of $\theta$ does every segment remain inside its
assigned square? The endpoint conditions alone do not ensure this:
for $\theta=9/10$ and $\mu=2$,
\[
    g_{9/10,2}(1/2) \ = \ -\frac{41}{800} \ < \ 0.
\]
Thus a connector can leave the square even though its endpoint
positions and derivatives are correct; see
Figure~\ref{fig:spiral-containment}.\\

\item \textbf{Self-intersection.}
Theorem~\ref{thm:spiral-global} rules out self-intersection
at the midpoint splice. For which other values of $\theta$ is
$\Gamma_\theta(\omega)$ simple for every placement word? If a
connector leaves its square, the containment argument no longer
applies, but this alone does not prove that the curve intersects
itself. If intersections can occur, what is their probability under
independent fair placements?\\

\item \textbf{Curvature propagation.}
Instead of prescribing entrance curvature from the scale $F_{n-1}$,
one could try to match the actual terminal curvature of the preceding
segment. Which corner configurations permit this while preserving
containment and uniform length bounds? Retaining the usual circular
arcs still prevents global $C^2$ continuity, since their curvatures
jump when the radii change.\\

\item \textbf{Asymptotic shape.}
How do winding, radial growth, and rescaled shape depend on the
placement word? What limiting behavior occurs under independent
fair placements, beyond the uniform length bounds proved above?\\

\item \textbf{Online constructions.}
The present corner rule needs one-step lookahead because the exit
corner and tangent depend on the next placement. Can a different
endpoint or tangent rule produce a contained, simple curve whose
segment is finalized as soon as its square is placed?\\

\item \textbf{Higher dimensions.}
Is there a natural curve, surface, or family of splines associated
with the order-$d$ tilings that reduces to this construction in
the plane?
\end{enumerate}

\section{Conclusion}\label{sec:conclusion}

We introduced a random growth model in which Fibonacci-sized blocks are attached along cyclically chosen coordinate axes while independent coin flips select the direction of each attachment.  The dimensions of the growing block are deterministic, but its location is governed by weighted Bernoulli subset sums.  Splitting the order-$d$ Fibonacci sequence into its coordinate residue classes exposes the main structure of the model: each class is superincreasing.  This gives unique subset sums, a most-significant-bit order, and canonical threshold words that convert geometric covering events into exact binary counts.

This representation leads to exact finite-time laws for the coordinate covering times and, by independence and inclusion--exclusion, for the full $d$-dimensional covering time.  For fixed $d\geq2$ and $I(a)\neq\varnothing$, the dominant root $\rho_d$ determines the macroscopic scale through
\begin{equation}
\label{eq:conclusion-covering-scale}
\mathbb{E}[T]\ =\ \log_{\rho_d}\delta_\infty(a)+O_d(1).
\end{equation}
At the same time, the centered fluctuations have exponential tails and moments bounded uniformly over the target.  Thus the covering time becomes relatively concentrated even though bounded lattice and phase effects may prevent a single centered limiting distribution without passing to subsequences.

The random-spiral construction gives a geometric companion to the covering theory.  Standard transitions retain their circular arcs, while nonstandard transitions are repaired by uniquely determined quintic connectors.  This produces a globally $C^1$ curve and recovers the usual Fibonacci spiral on the standard tiling word.  The construction should be viewed as a starting point rather than a final geometric model: the splice parameter remains free, and containment, simplicity, winding, and optimality of the resulting curves are not yet resolved.

Several probabilistic questions also remain.  It can be an exercise for the reader to find an explicit formula for when we hit the integer point $(i,j)$
in the regular, deterministic construction where the squares are added in the sequential pattern right, up, left, down. It would be natural to classify the subsequential laws of the centered covering time, to allow biased or dependent directional choices, and to determine which fixed-$d$ conclusions remain uniform when the dimension grows with the target.  Geometric questions are outlined in Section~\ref{section:questions}.

\section*{Acknowledgements}
This research was supported with funding from the National Science Foundation (grant DMS2341623), University of Rochester, Williams College Finnerty Fund, Missouri State University, and the University of Michigan. We would like to thank L\'aslz\'o Szalay for his ideas regarding the construction of the random Fibonacci spiral and Mateo Palomares for his contribution to the random spiral image. AI has been used in this project for proof-reading, assisting with LaTex formatting, and ensuring rigor.



\section*{Disclosure statement}
No conflict of interest has been reported by the authors.

\printbibliography

\medskip
\noindent MSC2020: 11B39.

\appendix
\section{Notation and conventions}\label{app:notation}

This appendix collects the notation used throughout the paper.  Unless a
different range is stated, $d\geq2$, coordinate indices satisfy
$0\leq r<d$, and all move and update indices are integers.  For $m_r,\tau_r,G_r,b_r$ and $B(a)$, assume $I(a)\neq\varnothing$ and use only active coordinates. The symbol $n$
usually denotes the number of completed global moves, $j$ an update number
within one coordinate (beginning with $j=0$), $m$ a number of available
coordinate weights, and $k$ either an auxiliary integer or a delay measured
in complete coordinate updates.

\subsection{General conventions}

\begin{description}
\item[$\mathbb N,\mathbb Z,\mathbb Q,\mathbb R,\mathbb C$]: 
The positive integers, integers, rational numbers, real numbers, and complex
numbers, respectively.  Indices that may equal zero are stated explicitly.

\item[$\mathbb P,\mathbb E,\operatorname{Var}$]:
Probability, expectation, and variance.

\item[$\#S$ and $|S|$]:
The cardinality of a finite set $S$.  Absolute values and complex moduli are
also written $|\cdot|$, with the meaning determined by context.

\item[$\mathbf 1_E$]:
The indicator of the event or condition $E$.

\item[$O_d(1),O_{d,q}(1),O_\theta(\cdot)$]:
The implied constant may depend only on the displayed subscript(s), and not
on the target point or placement word.  The symbols $o(1)$ and $\sim$ have their standard asymptotic meanings. The notation $\Theta_\theta$ means upper and lower bounds by positive constants depending only on $\theta$.

\item[$L^q$]:
The usual space of random variables with finite $q\textsuperscript{th}$ absolute moment; $q$
is a fixed real number with $q\geq1$ in the moment results.
\end{description}

\subsection{The tiling and the order-\texorpdfstring{$d$}{d} recurrence}

\begin{description}
\item[$F_n$]:
The ordinary Fibonacci sequence in the convention
$F_0=F_1=1$ and $F_n=F_{n-1}+F_{n-2}$ for $n\geq2$.

\item[$F_n^{(d)}$]:
The order-$d$ sequence defined by
$F_0^{(d)}=\cdots=F_{d-1}^{(d)}=1$ and
$F_n^{(d)}=F_{n-1}^{(d)}+F_{n-d}^{(d)}$ for $n\geq d$.

\item[$\mathcal R_n^{(d)}$]:
The random bounding hyperrectangle after $n$ completed moves; in particular,
$\mathcal R_0^{(d)}=[0,1]^d$.

\item[$c_n$]:
The independent fair placement bit used at global move $n$:
$0$ selects the positive face (right/up in the plane) and $1$ selects the
negative face (left/down in the plane).

\item[$\chi_d(x)$]:
The characteristic polynomial
$\chi_d(x)=x^d-x^{d-1}-1$.

\item[$\rho_d$]:
The unique positive root of $\chi_d$ and its unique root of largest modulus.
Equivalently, $\rho_d^{d-1}(\rho_d-1)=1$.

\item[$\varphi$]:
The golden ratio, $\varphi=\rho_2=(1+\sqrt5)/2$.

\item[$\zeta_1,\ldots,\zeta_d$]:
The roots of $\chi_d$, labeled so that $\zeta_d=\rho_d$.

\item[$\alpha_d$]:
The largest modulus of a nondominant root:
$\alpha_d=\max\{|\zeta|:\chi_d(\zeta)=0,\ \zeta\ne\rho_d\}$.

\item[$D_d(z)$]:
The generating-function denominator $D_d(z)=1-z-z^d$.

\item[$A_1,\ldots,A_d$]:
The coefficients in
$D_d(z)^{-1}=\sum_{j=1}^d A_j(1-\zeta_jz)^{-1}$; in particular,
$A_d=\rho_d^d/(\rho_d^{d-1}+d)$ is the dominant coefficient.

\item[$z_d$]:
The dominant pole of the generating function, $z_d=\rho_d^{-1}$.

\item[$n_d(x)$]:
The first recurrence index at which the sequence reaches $x$:
$n_d(x)=\min\{n\geq0:F_n^{(d)}\geq x\}$.
\end{description}

\subsection{Targets, coordinate weights, and canonical words}

\begin{description}
\item[$a=(a_0,\ldots,a_{d-1})$]:
The fixed lattice target in $\mathbb Z^d$.

\item[$\delta_r(a)$]:
The distance of the target beyond the relevant face of the initial cube in
coordinate $r$:
$\delta_r(a)=\max\{a_r-1,-a_r,0\}$.

\item[$I(a)$]:
The active-coordinate set
$I(a)=\{r:\delta_r(a)>0\}$.  When the target is fixed, the paper abbreviates
$\delta_r=\delta_r(a)$ and $I=I(a)$.

\item[$\delta_\infty(a)$]:
The largest coordinate distance,
$\delta_\infty(a)=\max_{0\leq r<d}\delta_r(a)$.

\item[$T(a)$]:
The first completed-move count for which $a\in\mathcal R_n^{(d)}$:
$T(a)=\inf\{n\geq0:a\in\mathcal R_n^{(d)}\}$.

\item[$w_{r,j}$]:
The thickness appended at the $j\textsuperscript{th}$ update of coordinate $r$ is
\begin{equation}
\label{eq:appendix-coordinate-weight}
w_{r,j}\ =\ F_{d-1+r+jd}^{(d)}.
\end{equation}
That update occurs at global move $r+1+jd$.

\item[$\mathcal W_{r,m}$ and $W_{r,m}$]:
The first $m$ weights in coordinate $r$ and their sum:
$\mathcal W_{r,m}=\{w_{r,0},\ldots,w_{r,m-1}\}$ and
$W_{r,m}=\sum_{j=0}^{m-1}w_{r,j}$, with
$\mathcal W_{r,0}=\varnothing$ and $W_{r,0}=0$.

\item[$u_r(n)$]:
The number of updates of coordinate $r$ among the first $n$ global moves is
\begin{equation}
\label{eq:appendix-update-count}
u_r(n)\ =\ \#\{0\ \leq\ i\ <\ n:i\equiv r\pmod d\}.
\end{equation}
Equivalently,
\begin{equation}
\label{eq:appendix-update-count-cases}
u_r(n)\ =\ \begin{cases}
0,&n\ \leq\ r,\\
1+\left\lfloor (n-1-r)/d\right\rfloor,&n\ \geq\ r+1.
\end{cases}
\end{equation}

\item[$\xi_{r,j}$]:
The target-directed indicator for the $j\textsuperscript{th}$ update of coordinate $r$: it is
$1$ precisely when that update is directed toward the target.  According to
the side on which the target lies, it is either $c_{r+1+jd}$ or
$1-c_{r+1+jd}$ from the global placement sequence.

\item[$E_{r,n}$] : 
The random subset of the weights available after $n$ moves that were directed
toward the target:
$E_{r,n}=\{w_{r,j}:0\leq j<u_r(n),\ \xi_{r,j}=1\}$.

\item[$x(E)$ and $x_j$]:
For $E\subseteq\mathcal W_{r,m}$, its membership word
$x(E)=(x_{m-1},\ldots,x_0)\in\{0,1\}^m$, where
$x_j=\mathbf1_{\{w_{r,j}\in E\}}$.

\item[$V_r(x)$]:
The weighted value of a binary word:
$V_r(x)=\sum_{j=0}^{m-1}x_jw_{r,j}$.

\item[$k(x,y)$ and $<_{\mathrm{MSB}}$]:
For distinct binary words, $k(x,y)=\max\{j:x_j\ne y_j\}$ is their largest
differing digit, and $x<_{\mathrm{MSB}}y$ means
$x_{k(x,y)}<y_{k(x,y)}$.

\item[$\operatorname{rank}_2(x)$]:
The ordinary binary rank
$\operatorname{rank}_2(x)=\sum_{j=0}^{m-1}x_j2^j$.

\item[$m_r=m_r(\delta_r)$]:
The first feasible coordinate-update count:
$m_r=\min\{m\geq1:W_{r,m}\geq\delta_r\}$.

\item[$x_r^\star$ and $x_{r,j}^\star$]:
The first word in MSB order whose weighted value reaches the target,
$x_r^\star=\min_{<_{\mathrm{MSB}}}\{x\in\{0,1\}^{m_r}:
V_r(x)\geq\delta_r\}$, and its digits.

\item[$R_r$]:
The binary rank of the canonical word,
$R_r=\operatorname{rank}_2(x_r^\star)$.

\item[$\kappa_r$]:
The normalized threshold rank
$\kappa_r=R_r/2^{m_r}\in[1/2,1)$.

\item[$v_{r,k}$]:
In the greedy construction, the contribution from the already fixed digits
larger than $k$:
$v_{r,k}=\sum_{j=k+1}^{m_r-1}x_{r,j}^\star w_{r,j}$.

\item[$\tau_r$]:
The number of coordinate-$r$ updates required for coverage, with $\inf\varnothing=\infty$:
$\tau_r=\inf\{m\geq0:\sum_{j=0}^{m-1}\xi_{r,j}w_{r,j}\geq\delta_r\}$.

\item[$G_r$]:
The additional coordinate-update delay after feasibility begins,
$G_r=\tau_r-m_r$.

\item[$b_r$]:
The earliest global move at which coordinate $r$ can be covered,
$b_r=r+1+d(m_r-1)$.

\item[$T_r$]:
The first global move at which active coordinate $r$ is covered:
$T_r=r+1+d(\tau_r-1)=b_r+dG_r$. For $r\notin I(a)$, set $T_r=0$.

\end{description}

\subsection{The full covering law and its moments}

\begin{description}
\item[$\mathcal G_r$]:
The sigma-algebra generated by the placement bits for coordinate $r$; see \eqref{eq:coordinate-sigma-algebra}.

\item[$\overline F_r(n)$]:
The coordinate survival function,
$\overline F_r(n)=\mathbb P(T_r>n)$.

\item[$T$]:
The full covering time when the target is fixed:
$T=T(a)=\max_{r\in I}T_r$, with the maximum over the empty set defined to be
$0$.

\item[$m_*$]:
The largest active threshold length,
$m_*=\max_{r\in I}m_r$.

\item[$n=dk+s$]:
The decomposition of global time into complete cycles and phase, where
$k\geq0$ and $0\leq s<d$.

\item[$J$]:
A nonempty subset of the active-coordinate set $I$ in the
inclusion--exclusion formulas.

\item[$R_J$]:
The product $R_J=\prod_{r\in J}R_r$.

\item[$\eta_s(J)$]:
The number of coordinates in $J$ updated before phase $s$:
$\eta_s(J)=\#\{r\in J:r<s\}$.

\item[$H_J$ and $K_J$]:
The phase sums
\begin{equation}
\label{eq:appendix-phase-sums}
H_J\ =\ \sum_{s=0}^{d-1}2^{-\eta_s(J)},
\qquad
K_J\ =\ \sum_{s=0}^{d-1}(2s+1)2^{-\eta_s(J)}.
\end{equation}

\item[$E_{\mathrm{fin}}$ and $E_{\mathrm{tail}}$]:
The finite and geometric-tail portions of the tail sum for $\mathbb E[T]$,
split at $dm_*$.

\item[$M_{2,\mathrm{fin}}$]:
The finite transient contribution to the second-moment tail sum.

\item[$z_J$]:
The geometric ratio $z_J=2^{-|J|}$.
\end{description}

\subsection{Asymptotic covering notation}

\begin{description}
\item[$B(a)$]:
The deterministic earliest feasible full-covering time,
$B(a)=\max_{r\in I(a)}b_r$.

\item[$G_{\max}$]:
The largest active coordinate delay,
$G_{\max}=\max_{r\in I(a)}G_r$.

\item[$\operatorname W_0$]:
The principal real branch of the Lambert function, characterized for $x>0$
by
\begin{equation}
\label{eq:appendix-lambert-function}
\operatorname W_0(x)e^{\operatorname W_0(x)}\ =\ x.
\end{equation}

\item[$m,u,v$ in the large-$d$ proof]:
Proof-local variables: $m=d-1$, $u=m(\rho_d-1)$, and $v$ is the argument of
the monotone function $v\mapsto v+\log v$.
\end{description}

\subsection{Random-spiral notation}

\begin{description}
\item[$\mathcal Q_n$]:
The square of side length $F_n$ in the planar tiling.

\item[$\omega=(\omega_1,\omega_2,\ldots)$]:
The binary placement word for the squares in the spiral construction.

\item[$\gamma_n$]:
The curve segment assigned to $\mathcal Q_n$; $\gamma_0$ is the initial
quarter-circle.

\item[$\theta$]:
The splice parameter in $(0,1)$ for a nonstandard segment.

\item[$A_n$ and $B_n$]:
The incoming and current-square scales in a nonstandard square:
$A_n=F_{n-1}$ and $B_n=F_n$.

\item[$h_n$ and $\varrho_n$]:
The connector width and outgoing-arc radius,
$h_n=\theta B_n$ and $\varrho_n=(1-\theta)B_n$.

\item[$y_{n,-}$]:
The reference incoming arc
$y_{n,-}(x)=A_n-\sqrt{A_n^2-x^2}$, used only to prescribe the incoming
endpoint jet.

\item[$y_{n,+}$]:
The outgoing circular arc
$y_{n,+}(x)=\sqrt{\varrho_n^2-(x-h_n)^2}$.

\item[$P_n$]:
The quintic connector on $[0,h_n]$ satisfying the six prescribed endpoint
value, slope, and second-derivative conditions.

\item[$t$]:
The normalized connector coordinate $t=x/h_n$.

\item[$\sigma_5(t)$]:
The quintic smoothstep polynomial
$\sigma_5(t)=10t^3-15t^4+6t^5$.

\item[$\Gamma_\theta(\omega)$]:
The full random Fibonacci spiral
$\Gamma_\theta(\omega)=\bigcup_{n\geq0}\gamma_n$.

\item[$\mu_n$]:
The consecutive-square ratio $\mu_n=F_n/F_{n-1}$.

\item[$g_{\theta,\mu}(t)$]:
The normalized graph of the quintic connector.

\item[$\beta_n(t)$]:
The connector parametrization
$\beta_n(t)=B_n(\theta t,g_{\theta,\mu_n}(t))$.

\item[$\mathcal L_n^{\mathrm{con}}(\theta)$]:
The length of the quintic portion of a nonstandard segment.

\item[$\ell_n$]:
The full length $\ell_n=\operatorname{len}(\gamma_n)$ of the segment assigned
to $\mathcal Q_n$.

\item[$\Lambda(\theta,\mu)$]:
The dimensionless connector-length integral
$\int_0^1\sqrt{\theta^2+(g'_{\theta,\mu}(t))^2}\,dt$.

\item[$A_\theta,B_\theta$]:
Positive constants, depending only on $\theta$, for which
$A_\theta F_n\leq\ell_n\leq B_\theta F_n$.

\item[$\mathcal K(s)$]:
Curvature as a function of arclength $s$ in the open geometric questions.

\item[$G^1$]:
Geometric continuity of position and oriented unit tangent; speed is not prescribed.

\item[$\gamma_n^{\mathrm{ref}}$]:
The reference composite segment in \eqref{eq:nonstandard-segment}, before transport by a rigid motion.

\item[$C^1,C^2$]:
The standard differentiability classes.  In the manuscript these refer to
regularity at joins after a compatible parametrization has been specified.
\end{description}

\paragraph{Scope-sensitive letters.}
The similarly lettered quantities $A_d$, $A_n$, and $A_\theta$ are,
respectively, a recurrence coefficient, an incoming spiral scale, and a
length-bound constant.  Likewise, $B(a)$, $B_n$, and $B_\theta$ are the
earliest feasible covering time, the current square size, and a length-bound
constant.  They occur in separate arguments but are not mathematically
related.  The global placement bits $c_n$ and target-directed bits
$\xi_{r,j}$ are also distinct: the latter may be the former or its complement.

\end{document}